\documentclass[11pt] {article} 
  \usepackage{amsmath} 
  \usepackage{mathrsfs}
  \usepackage[backref]{hyperref}

\usepackage{pdfsync}

\newtheorem{theorem}{Theorem}[section]
\newtheorem{lemma}{Lemma}[section]
\newtheorem{corollary}{Corollary}[section]

\newtheorem{remark}{Remark}[section]
\newcommand{\eqnsection}{
   \renewcommand{\theequation}{\thesection.\arabic{equation}}
   \makeatletter
   \csname @addtoreset\endcsname{equation}{section}
   \makeatother}
 
\def \ov{\overline}

\def \be{\begin{equation}} 
\def \ee{\end{equation}}
\def \bt{\begin{theorem}} 
\def \et{\end{theorem}}
\def \bl{\begin{lemma}} 
\def \el{\end{lemma}} 
\def \bea{\begin{eqnarray}}
\def \eea{\end{eqnarray}}
\def \bas{\begin{eqnarray*}}
\def \eas{\end{eqnarray*}}

\def \al{\alpha} 
\def \bb{\beta}
 
\def \Ga{\Gamma}

\def \De{\Delta} 
\def \ep{\epsilon}

\def \la{\lambda}

\def \om{\omega}
\def \Om{\Omega}

\def \si{\sigma}

\def \th{\theta}

\def \ze{\zeta}

\def \ff{\infty}
\def \wh{\widehat}
\def \wt{\widetilde}

\def \rar{\rightarrow}

\def \FF{{\cal F}}
\def \GG{{\cal G}}

\def \LL{{\cal L}}

\def \PP{{\cal P}}

\def \({\left(}
\def \){\right)}

\def \lc{\left\{}
\def \rc{\right\}}

\def \nn{\nonumber}
\def \Proof{\noindent{\bf Proof $\,$ }}

\def \bc{\begin{center} }
\def \ec{\end{center} }
\def \bs{\begin{slide} }
\def \es{\end{slide} }

\def\square{{\vcenter{\vbox{\hrule height.3pt
        \hbox{\vrule width.3pt height5pt \kern5pt
           \vrule width.3pt}
        \hrule height.3pt}}}}
\def\qed{{\hfill $\square$ \bigskip}}

\eqnsection
\begin{document}

\title{Markovian  loop soups: permanental processes and isomorphism theorems}

    \author{P.J. Fitzsimmons \,\,\, Jay Rosen \thanks{Research of the second author was  supported by  grants from the National Science Foundation, PSCCUNY and  grant number 208494  from the Simons Foundation.}}
\maketitle

\begin{abstract}   
   We construct  loop soups for general Markov processes without transition densities and show that the associated permanental process is equal in distribution to the loop soup local time. This is used to establish isomorphism theorems connecting the local time of the original process with the associated permanental process. Further properties of the loop measure are studied.   
 \end{abstract}

\bibliographystyle{amsplain}

  \section{Introduction}   \label{sec-1}
  
  A Markovian loop soup is a particular Poisson point process $\LL$ on paths associated to a Markov process $X$. It is determined by its intensity measure $\mu$ which we refer to as   loop measure.    Loop measure for Brownian motion was introduced by Symanzik in his seminal  paper \cite{Sy} on Euclidean quantum field theory, where it is referred to as `blob measure', and is a basic building block in his construction of quantum fields.   Brownian loop soup was introduced by Lawler and   Werner  \cite{LW},  in their work on SLE and conformally invariant processes in the plane. Le Jan extended the notion of loop soups to other Markov processes \cite{Le Jan}, and this has been generalized  further in \cite{LMR1, LMR2}. In all this work the loop measure is constructed using bridge measures for $X$. This requires that $X$ have transition densities. The main point of this paper is to show how to construct loop measures and hence loop soups for Markov processes
which have potential densities but not transition densities.

Our motivation in studying Markovian  loop soups is to better understand the wonderful and mysterious Isomorphism Theorem of Dynkin,  \cite{D83,D84}, which connects the family of total  local times $L=\{L^{x}_{\ff}, x\in S\}$ of a symmetric  Markov process  $X$ in $S$ with the Gaussian process $G=\{G_{x}, x\in S\}$ of covariance $u(x,y)$. (When $X$ is symmetric, $u(x,y)$ is positive definite.) Actually, in the Isomorphism Theorem it is  the 
family of squares of $G$, that is $G^{2}=\{G^{2}_{x}, x\in S\}$,
 which is connected with $L$. This theorem is not an isomorphism in the usual sense, but the connection between $L $ and $G^{2}$ is sufficiently tight that it has been used to derive many new properties of the local times, as described in   \cite{book}. This is why we consider the Isomorphism Theorem to be wonderful. We call it mysterious because it is hard to see intuitively why there should be any connection between Markov local times and Gaussian processes. 

As noted by Le Jan,  \cite[Theorem 9]{Le Jan1}, loop soups offer a deep understanding of this connection. Recall that each realization of $\LL$ is a countable collection of paths $\om $. Set
\begin{equation}
\wh L^{x}=\sum_{\om\in\LL}L^{x}_{\ff}(\om).\label{int.1}
\end{equation}
We call $\wh L^{x}$ the loop soup local time at $x$.
A simple application of the Palm formula for Poisson point processes provides a connection, an 
Isomorphism Theorem, between $L=\{L^{x}_{\ff}, x\in S\}$ and $\wh L=\{\wh L^{x}, x\in S\}$. Since $\wh L$ is defined  in terms of local times of $X$ this should not be surprising. 
What may be surprising is that when $X$ is symmetric then  $\wh L=\{\wh L^{x}, x\in S\}$ has the distribution of $G^{2}=\{G^{2}_{x}, x\in S\}$!
Furthermore, the definition of (\ref{int.1}) of $\wh L$ does not require the symmetry of $X$, so we obtain an Isomorphism Theorem for non-symmetric $X$.

In  1997,      D. Vere-Jones,  \cite{VJ}, introduced 
the $\al$-permanental process $\th:=\{\th_{x}, x\in S\} $ with kernel $u(x,y)$,   which is a real valued positive stochastic process   with 
joint distributions that satisfy 
\be
E\( \prod_{j=1}^{n}\th_{x_{j}}   \)=  \sum_{\pi\in \mathcal{P}_{n}}\al^{c(\pi)} \prod_{j=1}^{n}u(x_{j},x_{\pi(j)})\label{ls.14},
\ee
for any $x^{1},\ldots, x_{n}\in S$,
where $c(\pi)$ is the number of cycles in the permutation $\pi$ of $[1,n]$.  In addition, by  \cite[p. 128]{VJ},  the joint moment generating function of $ ( \th_{x_{1}},\ldots,  \th_{x_{n}}) $ has a non-zero radius of convergence. Consequently,  an  $\al$-permanental process is determined by its moments. It is not hard to show that in the symmetric case  $G^{2}/2=\{G^{2}_{x}/2, x\in S\}$ is a $1/2$-permanental process  with kernel $u(x,y)$, the covariance of $G$. 

In \cite{EK}, Eisenbaum and Kaspi  were able to show the existence of an $\al$-permanental process  with kernel $u(x,y)$ whenever 
$u(x,y)$ is the potential density of a transient Markov process $X$, and use this to obtain  an Isomorphism Theorem for non-symmetric $X$, where the role played in the symmetric case by the Gaussian squares $G^{2}$ is now played by a permanental process.   In this paper we will see that the loop soup local time $\wh L$   is an $\al$-permanental process  with kernel $u(x,y)$. 

The advantage of using loop soups to construct permanental processes and obtain Isomorphism Theorems is two-fold. First, as mentioned, loop soups provide an intuitive understanding of the connection between permanental processes and local times.
Second, this approach is capable of great generalization.  Recent work, \cite{LMR1,LMR2}, uses loop soups for Markov processes with potential densities $u(x,y)$ which may be infinite on the diagonal. In this case there are no local times and no permanental processes. Rather,  loop soups are used to prove the existence of permanental fields (indexed by measures rather than points in $S$)   
with which  
 to establish  Isomorphism Theorems: for continuous additive functionals in \cite{LMR1}, and for intersection local times in \cite{LMR2}. We know of no  way other than using loop soups  to  prove the existence of permanental fields associated with not necessarily symmetric  $X$, and the Isomorphism Theorems  contain constructs which seem inaccessible without the   loop soup. For example, in the symmetric case the Isomorphism Theorems  contain random variables  which are not in the associated Gaussian sigma field.  

Here is an outline of this paper.
The loop measure   is constructed and studied in Section \ref{sec-measure}. In the short sub-section \ref{sec-rot} we show that when transition densities exist,  our definition of loop measure agrees with the standard definition using bridge measures.
In Section \ref{sec-soup} we introduce the loop soup and quickly  show that the loop soup local time $\hat L$ is an $\al$-permanental process  with kernel $u(x,y)$. In the short  Section \ref{sec-LIT} we use the Palm formula to prove our Isomorphism Theorem. Further properties of the loop measure are derived in Sections \ref{sec-inv}-\ref{sec-trans}. These include invariance under loop rotation, and the behavior 
of the loop measure under restriction and space-time transformations. Here again the novelty is in deriving these properties in great generality and  without the assumption of transition densities.

Acknowledgements: We would like to thank Yves Le Jan for pointing out the connection between loop soups and isomorphism theorems.

\section{The loop measure}\label{sec-measure}

 Let $S$ be a locally compact 
 topological space
 with a countable base. 
  Let    $X\!=\!
(\Om,  \FF_{t}, X_t,\th_{t} ,P^x
)$ be a transient Borel right process with state space $S$,  and continuous potential densities 
$u(x,y)$ with respect to some $\si$-finite measure $m$ on $S$. That is 
\begin{equation}
P^x\(\int_0^\infty f(X_t)\,dt\) =\int_S u(x,y) f(y)\,m(dy),\qquad \forall x\in S,\nn
\end{equation} 
for each non-negative Borel function $f:S\to[0,\infty)$. We assume further that $u(x,y)>0$ for all $x,y\in S$. This amounts to the assumption that each point is regular, and that the process is irreducible in the sense that $P^x(T_y < \infty)>0$ for all $x,y$. Then  there is a function $(\omega,s,y)\to L^y_s(\omega)$ from $\Omega\times[0,\infty)\times S$ to $[0,\infty)$ that is jointly progressively measurable in $(\omega,s)$ and Borel measurable in $y$, such that 
 for each $y\in S$, $t\mapsto L^y_t$ is a CAF increasing only when $X$ is in state $y$ ({\it i.e.}, a local time at $y$), and $E^{x}(L_{\ff}^{y})=u(x,y)$ for all $x,y\in S$. This follows from the proofs in \cite[Section 3.6]{book} if we choose the approximate delta functions $f_{\ep,y}(x)$ used there to be  of the form
 \begin{equation}
 f_{\ep,y}(x)={f_{\ep}(d(y,x)) \over \int_S f_{\ep}(d(y,z))\,m(dz)}\label{met1}
 \end{equation}
 where $d$ is a metric for the topology of $S$, and $f_{\ep}$ is a continuous function supported on $[0,\ep]$, and define $L^y_t(\om)=\liminf_{n\to\ff}\int_0^t f_{n^{-1},y}(X_s)\,ds$. (This is used to show measurability in $y$).

 Under our assumption that $u(x,y)$ is continuous, it follows as in the proof of \cite[Lemma 3.4.3]{book}, that uniformly in $x$, $u(x,y)$ as a function of $y$ is locally bounded and continuous. This implies that for any $\bb>0$, the same is true for
\begin{equation}
u^{\bb}(x,y):=u(x,y)-\bb\int U^{\bb}(x,dz)u(z,y),\label{res.1}
\end{equation}
and it follows from the resolvent equation that for each $x$, $u^{\bb}(x,y)$ is a   density for $U^{\bb}(x,dy)$ with respect to $m(dy)$. It then follows from the resolvent equation that for any $\al,\bb>0$ and all $x,y$
\begin{equation}
{u^{\al}(x,y)-u^{\bb}(x,y) \over \bb-\al}=\int u^{\al}(x,z)u^{\bb}(z,y)\,dm(z)=\int u^{\bb}(x,z)u^{\al}(z,y)\,dm(z).\label{res.2}
\end{equation}
Using (\ref{res.1}) and the resolvent equation for additive functionals we see that
\begin{equation}
u^{\bb}(x,y)=E^{x}\(\int_{0}^{\ff}e^{-\bb t}d_{t}L^{y}_{t}\).\label{res.3}
\end{equation}

We now show that  $m(K)<\ff$ for each compact $K\subseteq S$. To see this note first that from (\ref{res.3})  and our assumption that $u(x,y)>0$ for all $x,y\in S$, that also  $u^{1}(x,y)>0$ for all $x,y\in S$.   It then follows from the last paragraph that  $y\mapsto u^1(x,y)$ is bounded below  for $y\in K$
by a constant $C = C(x)>0$. Consequently 
\be
C\cdot m(K)\leq \int_K u^1(x,y)\,m(dy)  \leq \int_S u^1(x,y)\,m(dy)=U^11(x) \leq  1.\label{radon}
\ee

 We may take the canonical representation of $X$ in which $\Om$ is  the set of right continuous paths $\om$ in $S_{\De}=S\cup \De$ with  $\De\notin S$, and  is such that   $\om(t)=\De$ for all $t\geq \ze=\inf \{t>0\,|\om(t)=\De\}$.   Then $X_{t}(\om)=\om(t)$.   We define a $\si$-finite measure $\mu$
on $(\Om, \mathcal{F})$ by 
\begin{equation}
\int F\,d\mu=\int_{S}P^{x}\(\int_{0}^{\ff}{1\over t}\,F\circ k_{t}\, \,d_{t}L_{t}^{x}\)\,dm (x),  \label{ls.3}
\end{equation}
for all $ \mathcal{F}$-measurable functions $F$  on $\Om$.
Here $k_{t}$ is the killing operator defined by $k_{t}\om(s)=\om(s)$ if $s<t$ and $k_{t}\om(s)=\De$ if $s\geq t$. 
  We call $\mu$ the loop measure of $X$ because, when $X$ has continuous paths, $\mu$ is concentrated on the set of continuous loops. 
See also  Lemma \ref{lem-ls3} below. Even if $X$ is not assumed to have continuous paths
we can  verify that $\mu$ is concentrated on $\{X_{0}=X_{\ze^{-}}, \ze<\ff\}$. To see this, note first of all that since $1_{\{\ze=\ff\}}\circ k_{t}=0$ for each $t$, it is clear from (\ref{ls.3}) that $\mu(\ze=\ff)=0$. Then, since $L_{t}^{x}$ is constant for $t\geq \ze$, while on $t\leq \ze$ we have $\ze\circ k_{t}=t$,
\be 
 P^{x}\(\int_{0}^{\ff}{1\over t}\,1_{\{X_{0}\neq X_{\ze^{-}}\}}\circ k_{t}\, \,d_{t}L_{t}^{x}\)
 =P^{x}\(\int_{0}^{\ff}{1\over t}\,1_{\{x\neq X_{t^{-}}\}}\, \,d_{t}L_{t}^{x}\). \label{atom1}
\ee
But by right-continuity of paths, the set of times for which $X_{t-}(\omega) $ either fails to exist or exists but is different from $X_t(\omega)$ is at most countable, for each $\omega\in \Omega$, \cite[IV, Theorem 88D]{DM}, while $L_{t}^{x}$ is continuous in $t$ so that $d_{t}L_{t}^{x}$ has no atoms. Hence (\ref{atom1})
\begin{equation}
= P^{x}\(\int_{0}^{\ff}{1\over t}\,1_{\{x\neq X_{t}\}}\, \,d_{t}L_{t}^{x}\)=0,  \label{atom2}
\end{equation}
where the last equality used the fact that  $d_{t}L_{t}^{x}$ is supported on $\{ X_{t}=x\}$.

As usual, if $F$ is a function, we often write $\mu (F)$ for $\int F\,d\mu$.

\begin{lemma}\label{lem-ltmom} For any $k$, and any $y_{1},\ldots,y_{k}\in S$
\begin{equation}
\mu\(\prod_{j=1}^{k}L_{\ff}^{y_{j}}\)= \sum_{\pi\in \mathcal{P}_{k-1}}   u(y_{k},y_{\pi(1)})\cdots  u(y_{\pi(k-2)},y_{\pi(k-1)}) u(y_{\pi(k-1)},y_{k}),\label{ls.5c4}
\end{equation}
where $\mathcal{P}_{k-1}$ denotes the set of permutations of $[1,k-1]$.  When $k=1$ this means 
$\mu\(L_{\ff}^{y_{1}}\)=u(y_{1},y_{1})$.
 \end{lemma}
 
 \Proof We present a derivation of Lemma \ref{lem-ltmom} suggested by Symanzik, \cite{Sy}.
We first show  that for any $k,\al\geq 0$ and $x,y\in S$
 \bea
 &&\mathcal{V}:=P^{x}\left( \int_0^\infty  e^{-\al t}\, \prod_{j=1}^{k} L_{t}^{z_{j}}\,\,   d_tL^{y}_t\right)\label{6.25}\\
&=& \sum_{\pi\in\PP_{k}}
u^{\al}(x, z_{\pi (1)})u^{\al}(z_{\pi (1)},z_{\pi (2)})\cdots u^{\al}(z_{\pi (k-1)},z_{\pi (k) })u^{\al}(z_{\pi (k)},y).\nn
 \eea
 To see this, note that   
\begin{equation}
\mathcal{V}=\sum_{\pi\in \mathcal{P}_{k}} 
 P^{x}\left(\int_0^\infty  e^{-\al t}\, \int_{\{0<s_{1}<\cdots<s_{k}<t\}} \, \prod_{j=1}^{k}\,dL_{s_{j}}^{z_{\pi(j)}}\,\,   d_tL^y_t\right),\label{t1j}
\end{equation}
and for each  $\pi\in \mathcal{P}_{k}$ we have 
 \bea
 \mathcal{J}_{\pi}&:=&
 P^{x}\left(\int_0^\infty e^{-\al t}\, \int_{\{0<s_{1}<\cdots<s_{k}<t\}} \, \prod_{j=1}^{k}\,dL_{s_{j}}^{z_{\pi(j)}}\,\,   d_tL^y_t\right)\label{t2j}\\
  & =&  P^{x}\left(\int_{\{0<s_{1}<\cdots<s_{k}<\ff\}} \(\int_{s_{k}}^\infty   \, \,\,e^{-\al t}\,   d_tL^y_t  \)    \prod_{j=1}^{k}\,d_{s_{j}}L_{s_{j}}^{y_{\pi(j)}}\right)  \nonumber\\
 &=&    P^{x}\left(\int_{\{0<s_{1}<\cdots<s_{k}<\ff\}}e^{-\al s_{k}}  \(\int_{0}^\infty   \, \,\,e^{-\al t}\,   d_tL^y_t  \) \circ \th_{s_{k}}   \,\,  \prod_{j=1}^{k}\,d_{s_{j}}L_{s_{j}}^{z_{\pi(j)}}\right).  \nonumber 
 \eea
Using the Markov property, see for example Theorems 28.7 and 22.8 of \cite{S}, and (\ref{res.3}), we have 
\bea
 \mathcal{J}_{\pi}&=& P^{x}\left(\int_{\{0<s_{1}<\cdots<s_{k}<\ff\}}e^{-\al s_{k}} E^{X_{s_{k}}} \(\int_{0}^\infty   \, \,\,e^{-\al t}\,   d_tL^y_t  \)   \,\,  \prod_{j=1}^{k}\,d_{s_{j}}L_{s_{j}}^{z_{\pi(j)}}\right)\nn\\
 &=&P^{x}\left(\int_{\{0<s_{1}<\cdots<s_{k}<\ff\}}    \,e^{-\al s_{k}}\,  \prod_{j=1}^{k}\,dL_{s_{j}}^{z_{\pi(j)}}\right)\,u^{\al}(z_{\pi(k)},y).\label{t4j}
\eea
It then follows by induction that
\begin{equation}
 \mathcal{J}_{\pi}=u^{\al}(x, z_{\pi (1)})u^{\al}(z_{\pi (1)},z_{\pi (2)})\cdots u^{\al}(z_{\pi (k-1)},z_{\pi (k) })u^{\al}(z_{\pi (k)},y).\label{t5j}
\end{equation}
Using (\ref{t1j}) then proves (\ref{6.25}).

It follows from (\ref{6.25}) that 
 \bea
 &&\int_{S}P^{x}\left( \int_0^\infty  e^{-\al t}\, \prod_{j=1}^{k} L_{t}^{z_{j}}\,\,   d_tL^{x}_t\right)\,dm (x)\label{6.26}\\
&=& \sum_{\pi\in\PP_{k}}
\int_{S} u^{\al}(x, z_{\pi (1)})u^{\al}(z_{\pi (1)},z_{\pi (2)})\cdots u^{\al}(z_{\pi (k-1)},z_{\pi (k) })u^{\al}(z_{\pi (k)},x)\,dm (x).\nn
 \eea
By (\ref{res.3}) we have that $u^{\bb}(x,y)\uparrow u^{\al}(x,y)$ as $\bb\downarrow  \al$. Hence by (\ref{res.2}) and the monotone convergence theorem
 \begin{equation}
\int_{S}  u^{\al}(y, x)u^{\al}(x, z) \,dm (x)=-  {d \over d\al} u^{\al}(y, z).  \label{6.27}
\ee
Hence  the right hand side of (\ref{6.26})
 \begin{equation}
=-\sum_{\pi\in\PP_{k}}
u^{\al}(z_{\pi (1)},z_{\pi (2)})\cdots u^{\al}(z_{\pi (k-1)},z_{\pi (k) }){d \over d\al}u^{\al}(z_{\pi (k)},z_{\pi (1)}). \label{6.28}
 \end{equation}
 The sum is over  all permutations of the `labels' of the points $z_{1},\ldots, z_{k}$ which in this expression appear in a circle. By fixing $z_{k}$ and considering permutations  $\pi\in\PP_{k-1}$, we can rewrite (\ref{6.28}) as
 \begin{eqnarray}
&& \hspace{-.3 in}-\sum_{\pi\in\PP_{k-1}}\sum_{j=1}^{k}
u^{\al}(z_{k},z_{\pi (1)})\cdots {d \over d\al}u^{\al}(z_{\pi (j-1)},z_{\pi (j)})\cdots u^{\al}(z_{\pi (k-1)},z_{k})
 \label{6.29}\\
 && \hspace{-.3 in} =-\sum_{\pi\in\PP_{k-1}} {d \over d\al}\(
u^{\al}(z_{k},z_{\pi (1)})\cdots u^{\al}(z_{\pi (j-1)},z_{\pi (j)})\cdots u^{\al}(z_{\pi (k-1)},z_{k})\), \nonumber
 \end{eqnarray} 
 where for notational convenience  we set $\pi(0)=\pi(k)=k$ for $\pi\in\PP_{k-1}$ in the first line. The second line is the product rule for differentiation.
 
By (\ref{res.3})
  \begin{equation}
 \lim_{\al\rar\ff} \,\,u^{\al}(x, y)=0.\label{6.10}
 \end{equation}
Hence by what we have shown about the right hand side of (\ref{6.26})
\begin{eqnarray}
&&\int_{0}^{\ff}\( \sum_{\pi\in\PP_{k}}
\int_{S} u^{\al}(x, z_{\pi (1)})u^{\al}(z_{\pi (1)},z_{\pi (2)})\cdots u^{\al}(z_{\pi (k)},x)\,dm (x)\)\,d\al
\nn\\
&&=-\sum_{\pi\in\PP_{k-1}}\int_{0}^{\ff} {d \over d\al}\(
u^{\al}(z_{k},z_{\pi (1)})\cdots u^{\al}(z_{\pi (j-1)},z_{\pi (j)})\cdots u^{\al}(z_{\pi (k-1)},z_{k})\)\,d\al
\nn\\
&&= \sum_{\pi\in\PP_{k-1}} 
u(z_{k},z_{\pi (1)})\cdots u(z_{\pi (j-1)},z_{\pi (j)})\cdots u (z_{\pi (k-1)},z_{k}). \label{6.31}
\end{eqnarray}
Thus by (\ref{6.26})
 \bea
 &&\int_{0}^{\ff} \lc \int_{S}P^{x}\left( \int_0^\infty  e^{-\al t}\, \prod_{j=1}^{k} L_{t}^{z_{j}}\,\,   d_tL^{x}_t\right)\,dm (x)\rc\,d\al\label{6.36}\\
&=& \sum_{\pi\in\PP_{k-1}} 
u(z_{k},z_{\pi (1)})\cdots u(z_{\pi (j-1)},z_{\pi (j)})\cdots u(z_{\pi (k-1)},z_{k}),\nn
 \eea
and Lemma \ref{lem-ltmom} follows by applying Fubini's theorem to interchange the order of integration.\qed

 Let $Q^{x,y}$ denote the   measure defined
on $(\Om, \mathcal{F})$ by 
\begin{equation}
Q^{x,y}(F)= P^{x}\(\int_{0}^{\ff} \,F\circ k_{t}\, \,d_{t}L_{t}^{y}\),  \label{ls.3n}
\end{equation}
for all $ \mathcal{F}$ measurable functions $F$ on $\Om$. Since $\zeta \circ k_{t}=\ze\wedge t$,
it follows  that if $F_s\in b\mathcal{F}^{0}_s$ where $\mathcal{F}^{0}_s$ is the $\sigma$-algebra generated by $\{X_{r},\,0\leq r\leq s\}$ then  $(1_{\{\ze>s \}}\,\,F_s)\circ k_{t}=1_{\{\ze\wedge t>s \}}\,\,F_s$. Hence
\bea
Q^{x,y}( 1_{\{\ze>s \}}\,\,F_s)&= &P^{x}\(F_s \int_{s}^{\ff}  \,d_{t}L_{t}^{y}\)\label{nit4.4}\\
&= &P^{x}\(F_s \,\, L_{\ff}^{y}\circ \th_{s}\)=P^x( F_s \,\,u(X_s,y )). 
\nn
\eea
We remark that under the measures $P^{x/h}={1 \over u(x,y)}Q^{x,y}$, the paths of $X$ are conditioned to hit $y$  and die on their last exit from $y$. $P^{x/h}$ is the $h$-transform of $P^{x}$,
with $h(x)=u(x,y)/u(y,y)=P^{x}(T_{y}<\ff)$.
 
In the proof of the Isomorphism Theorem we will need the following Lemma.
\bl\label{lem-qmu}For any $k$, $x, x_{j}\in  S $, $j=1,\ldots,k$, and any bounded measurable function $H$ on $R^{k}$ we have 
\begin{equation}
\mu\(L^{x}_{\ff}\,H\(L^{x_{1}}_{\ff},\cdots, L^{x_{k}}_{\ff}\)\)=Q^{x,x}\(H\(L^{x_{1}}_{\ff},\cdots, L^{x_{k}}_{\ff}\)\).\label{qmu}
\end{equation}
\el

{\bf   Proof  of Lemma \ref{lem-qmu}: } Let $x,y,x_{j}\in  S $, $j=1,\ldots,k$.  Since $L_{\ff}^{x_{j}}\circ k_{t}=L_{t}^{x_{j}}$, it follows from (\ref{ls.3n}) that
\begin{equation}
Q^{x,y}\(   \prod_{j=1}^{k}L^{x_{j}}_{\ff}\)=P^{x}\(\int_{0}^{\ff} \,  \prod_{j=1}^{k}L^{x_{j}}_{t}\, \,d_{t}L_{t}^{y}\).\label{geny}
\end{equation}
It then  follows from (\ref{6.25})   that
 \bea
&& 
 Q^{x,y}\(   \prod_{j=1}^{k}L^{x_{j}}_{\ff}\)=\sum_{\pi\in \mathcal{P}_{k}} 
 u(x, x_{\pi (1)})u(x_{\pi (1)},x_{\pi (2)})\cdots \label{nit8.2q0}\\
&&  \hspace{1.3 in}\cdots  u(x_{\pi (k-1)},x_{\pi (k) })u(x_{\pi (k)},y). \nonumber
\eea
Comparing (\ref{nit8.2q0}) with $y=x$ and (\ref{ls.5c4}) we see that (\ref{qmu}) holds for all polynomial $H$. But it is easily seen from (\ref{ls.5c4}) and (\ref{nit8.2q0}) that the random variables $L^{z}_{\ff}$ are exponentially integrable both under $Q^{x,x}$ and $\mu\(L^{x}_{\ff}\,\cdot\)$, hence finite dimensional distributions are determined by their moments.\qed

Since $\zeta \circ k_{t}=\ze\wedge t$,  we note for future reference that
\begin{eqnarray}
\mu(F)&=&\int_{S}P^{x}\(\int_{0}^{\ff}{1\over t}\, F\circ k_{t}\, \,d_{t}L_{t}^{x}\)\,dm (x)
\label{6.41}\\
&=&\int_{S}P^{x}\(\int_{0}^{\ff} \,\( {F \over \ze} \)\circ k_{t}\, \,d_{t}L_{t}^{x}\)\,dm (x) \nonumber\\
&=&\int_{S}Q^{x.x}\( {F \over \ze} \)\,dm (x).   \nonumber
\end{eqnarray}

In the sequel we will use the fact that $(t,x)\mapsto L^x_t(\omega)$ is an occupation density with respect to $m$:
\begin{equation}
\int_0^t f(X_s)\,ds =\int_S f(x) L^x_t\,m(dx), \label{odf}
\end{equation} 
for all $t\geq 0$ and all  non-negative Borel functions $f $, almost surely. It suffices to prove this for $f\geq 0$ which are continuous and compactly supported. This case follows from the proof of \cite[Theorem 3.7.1]{book}, with one change. That theorem assumed the joint continuity of $L^x_t$ in order to show that the right hand side of (\ref{odf}), which we denote by $A_t$, is a CAF. But this can be seen directly. $A_t$ is obviously monotone increasing in $t$ and constant for $t\geq \ze$. Also, using (\ref{radon})
\begin{equation}
E^{y}\(A_\ff\)=\int_S u(y,x)f(x) \,m(dx)<\ff,\label{odf1}
\end{equation}
hence $A_\ff<\ff $ a.s.  Hence the a.s. continuity of $A_t$ follows from  the dominated convergence theorem after applying Fubini to the fact  for each $x\in S$, a.s. in $\om$,   $L^x_{t}(\om)$ is continuous  in $t$. Finally, fix $s,t>0$. We have   $L^x_{s+t}(\om)=L^x_{s}(\om)+L^x_{t}\circ\th_{s}(\om)$
for each $x\in S$, a.s. in $\om$. Hence by Fubini this holds  a.s. in $\om$ for a.e. $x\in S$. From the right hand side of (\ref{odf}) we then see that a.s. in $\om$, $A_{s+t}(\om)=A_{s}(\om)+A_{t}\circ\th_{s}(\om)$, which completes the proof that $A_t$ is a CAF, and hence the proof of (\ref{odf}).

\subsection{Transition densities}\label{sec-rot}

 For this subsection only, we assume that  
$P_{t}(x,dy)\ll dm(y)$ for each $t>0$ and $x\in S$; in other words,  $P_{t}(x,dy)$ has transition densities with respect to $m$.   Under this assumption we give an alternate description of the loop measure. This is the description found in the literature.  Using this description we  give a simple proof of the fact that the loop measure is invariant under loop rotation.   A proof of this fact without assuming transition densities is given in Section \ref{sec-inv}.        The material in this sub-section will not be used in the following sections of the paper.

  Under our assumption that  $P_{t}(x,dy)$ has transition densities with respect to $m$, it   follows from \cite{D80} that we can find 
 jointly measurable transition densities $p_{t}(x,y)$ with respect to $m$ which satisfy the Chapman-Kolmogorov equation
 \begin{equation}
\int p_{s}(x,y)p_{t}(y,z)\,dm(y)=p_{s+t}(x,z).\label{CK}
 \end{equation}
Assume that $p_{t}(x,y)<\ff$, for all $0<t<\ff$ and $ x,y\in S$.
It then  follows as in  \cite{FPY} that for all $0<t<\ff$ and $x,y\in S$,   there exists a finite measure $Q_{t}^{x,y}$ on $\mathcal{F}_{t^{-}}$, of total mass $p_{t}(x,y)$, such that
\begin{equation}
Q_{t}^{x,y}\(F_{s}\)=P^{x}\( F_{s}  \,p_{t-s}(X_{s},y) \),\label{10.1}
\end{equation}
for all $F_{s} \in \mathcal{F}_{s}$ with $s<t$.  In particular, for any   $0< t_{1}\leq \cdots \leq t_{k-1}\leq t_{k} < t$ and bounded  Borel measurable functions $f_{1},\ldots, f_{k}$ 
\begin{eqnarray}
&&Q_{t}^{x,y}\( \prod_{j=1}^{k}f_{j}(X_{t_{j}})\)
\label{br7}\\
&&  =\int_{S^{k}} p_{t_{1}}(x,y_{1})f_{1}(y_{1}) p_{t_{2}-t_{1}}(y_{1},y_{2})f_{2}(y_{2})\cdots  \nonumber\\
&&\hspace{.2 in}\cdots p_{t_{k}-t_{k-1}}(y_{k-1},y_{k}) f_{k}(y_{k})p_{t-t_{k}}(y_{k},y)\,dm(y_{1})
\cdots \,dm (y_{k}).\nn
\end{eqnarray}

\begin{lemma}\label{lem-eq}
\begin{equation}
\mu(F)=\int_{0}^{\ff}{1\over t}\int_{S}\,Q_{t}^{x,x}\(F\circ k_{t} \)\,dm (x)\,dt\label{eq1}
\end{equation}
for all $ \mathcal{F}$ measurable functions F on $\Om$.
\end{lemma}

{\bf  Proof of Lemma \ref{lem-eq}} Let us temporarily use the notation $\wt\mu (F)$ to denote the right hand side of (\ref{eq1}). It suffices to show that   $\mu (F)=\wt\mu (F)$ for all $F$ of the form 
$F=\prod_{j=1}^{k}f_{j}(X_{t_{j}})$, for all
    $0< t_{1}\leq \cdots \leq t_{k-1}\leq t_{k} < \ff$ and bounded  Borel measurable functions $f_{1},\ldots, f_{k}$ on $S\cup \De$ with $f_{j}(\De)=0$, $j=1,\ldots,k$. Note that this last condition implies that
 \begin{equation}
\prod_{j=1}^{k}f_{j}(X_{t_{j}}\circ k_{t})=1_{\{t_{k}<t\}}\prod_{j=1}^{k}f_{j}(X_{t_{j}}).\label{10.5}
\end{equation}

Using  (\ref{br7})
\begin{eqnarray}
\wt\mu (F)&=&\int_{0}^{\ff}{1\over t}\int_{S}\,Q_{t}^{x,x}\(\(   \prod_{j=1}^{k}f_{j}(X_{t_{j}})\)\circ k_{t} \)\,dm (x)\,dt
\label{br8}\\
&=& \int_{t_{k}}^{\ff}{1\over t}\int_{S}\,Q_{t}^{x,x}\(   \prod_{j=1}^{k}f_{j}(X_{t_{j}})  \)\,dm (x)\,dt  \nonumber\\
&=&  \int_{t_{k}}^{\ff}{1\over t}\int_{S^{k+1}} p_{t_{1}}(x,y_{1})f_{1}(y_{1}) p_{t_{2}-t_{1}}(y_{1},y_{2})f_{2}(y_{2})\cdots p_{t_{k}-t_{k-1}}(y_{k-1},y_{k}) \nonumber\\
&&\hspace{.05 in}\cdots  f_{k}(y_{k})\,\,p_{t-t_{k}}(y_{k},x)\,dm(y_{1})
\cdots \,dm (y_{k})\,dm (x)\,dt.\nn
\end{eqnarray}

Similarly, using the Markov property
\bea
\mu\(F\)&=&\int_{S}P^{x}\(\int_{0}^{\ff}{1\over t}\,\prod_{j=1}^{k}f_{j}(X_{t_{j}})\circ k_{t} \,dL_{t}^{x}\)\,dm (x) \label{ls.1}\\
&=&\int_{S}P^{x}\(\int_{t_{k}}^{\ff}{1\over t}\,\prod_{j=1}^{k}f_{j}(X_{t_{j}})\,dL_{t}^{x}\)\,dm (x) \nn\\
&=&\int_{S}P^{x}\(\prod_{j=1}^{k}f_{j}(X_{t_{j}})\,\(\int_{t_{k}}^{\ff}{1\over t}\,d_{t}L_{t-t_{k}}^{x}\)\circ \th_{t_{k}}\)\,dm (x) \nn\\
&=&\int_{S}P^{x}\(\prod_{j=1}^{k}f_{j}(X_{t_{j}})\,E^{X_{t_{k}}}\(\int_{t_{k}}^{\ff}{1\over t}\,d_{t}L_{t-t_{k}}^{x}\)\)\,dm (x) \nn\\
&=& \int_{S^{k+1}} p_{t_{1}}(x,y_{1})f_{1}(y_{1}) p_{t_{2}-t_{1}}(y_{1},y_{2})f_{2}(y_{2})\cdots  p_{t_{k}-t_{k-1}}(y_{k-1},y_{k})\nonumber\\
&&\hspace{.2 in}\cdots  f_{k}(y_{k})\,E^{y_{k}}\(\int_{t_{k}}^{\ff}{1\over t}\,d_{t}L_{t-t_{k}}^{x}\)\,dm (x)\,dm(y_{1})
\cdots \,dm (y_{k}).\nn
\eea 

 For $y,x\in S$,  we define the measure $\Gamma_{y,x}(\cdot )$  on $[0,\infty)$ with cdf:
\be
\Gamma_{y,x}([0,t]):=E^y(L^x_t),\label{gamdef}
\ee
so that for all bounded measurable functions $g$ on $[0,\infty)$
\begin{equation}
\int_{0}^{\ff}g(t)\Gamma_{y,x}(dt)= E^y\(\int_{0}^{\ff}g(t)\,d_{t}L^x_t\).\label{intp}
\end{equation}

We claim that for each $y$, we can find a set $S_{y}\subset S$ with $m(S_{y})=0$ such that 
\begin{equation}
E^y(L^x_t)=\int_{0}^{t}p_{s}(y,x)\,ds\label{cla.1}
\end{equation}
for all $t$ and $x\in S^{c}_{y}$. 
Then by (\ref{intp}), for any  $t_{k}$   and all $x\in S^{c}_{y}$
 \begin{equation}
\int_{t_{k}}^{\ff}{1\over t}\,p_{t-t_{k}}(y_{k},x) \,dt=E^{y_{k}}\(\int_{t_{k}}^{\ff}{1\over t}\,d_{t}L_{t-t_{k}}^{x}\).\label{cla.2}
\end{equation}
Since the right hand side of (\ref{ls.1}) involves a $dm (x)$ integral, by Fubini we can replace the term $E^{y_{k}}\(\int_{t_{k}}^{\ff}{1\over t}\,d_{t}L_{t-t_{k}}^{x}\)$ which appears there with the left hand side of (\ref{cla.2}). Thus we will obtain
\begin{eqnarray}
 \mu (F)&=&  \int_{t_{k}}^{\ff}{1\over t}\int_{S^{k+1}} p_{t_{1}}(x,y_{1})f_{1}(y_{1}) p_{t_{2}-t_{1}}(y_{1},y_{2})f_{2}(y_{2})\cdots  \label{cla.3}\\
&&\hspace{.05 in}\cdots p_{t_{k}-t_{k-1}}(y_{k-1},y_{k}) f_{k}(y_{k})p_{t-t_{k}}(y_{k},x)\,dm(y_{1})
\cdots \,dm (y_{k})\,dm (x)\,dt.\nn
\end{eqnarray}
Comparing with (\ref{br8})  then shows that $\mu (F)=\wt\mu (F)$. 

It only remains to verify our claim concerning (\ref{cla.1}). Note that since the left hand  side of (\ref{cla.1}) is continuous in $t$ and the right hand side is monotone,  it suffices to find a set $S_{y}$  which works for all rational $t$, hence for each fixed $t$. By the occupation density formula (\ref{odf})
\bea
\int_S f(x)\( \int_{0}^{t}p_{s}(y,x)\,ds\)\,m(dx)&=&E^y\(\int_0^t f(X_s)\,ds \)\label{cla.5}\\
&=&\int_S f(x) \(E^y(L^x_t)\)\,m(dx).\nn
\eea
Since this holds for all bounded measurable $f$, our claim for fixed $t$ is established.
\qed

For later use we note that applying the Chapman-Kolmogorov equation (\ref{CK}) for the $dm (x)$ integral in (\ref{cla.3}) shows that 
\bea
\lefteqn{
\mu\(   \prod_{j=1}^{k}f_{j}(X_{t_{j}})\) \label{ls.4}}\\
&& =\int_{t_{k}}^{\ff}{1 \over t} \int_{S^{k}} f_{1}(y_{1}) p_{t_{2}-t_{1}}(y_{1},y_{2})f_{2}(y_{2})\cdots \nonumber\\
&&\hspace{.2 in}\cdots  p_{t_{k}-t_{k-1}}(y_{k-1},y_{k})f_{k}(y_{k})p_{t_{1}+t-t_{k}}(y_{k},y_{1})\,dm(y_{1})
\cdots \,dm (y_{k})\,dt.\nn
\eea

 The next result  justifies our calling  $\mu$ the loop measure even for  a process with discontinuous paths.  This result will be proved in full generality in section 5.
 Define the loop rotation $\rho_{u}$
by
\begin{equation}
\rho_{u}\om (s)
	=\left\{\begin{array}{ll}
	\om (s+u\mod\ze(\om)), & \mbox{ if $0\leq s<\ze(\om)$}\\
	\De, & \mbox{otherwise.}
	\end{array} \right.
	\label{rosen:d1.4c}
\end{equation}
Here, for two positive numbers $a,b$ we define  $a\mod b= a-mb$ for the unique positive integer $m$ such that $\,0\leq a-mb<b $ . Set $(a)_{b}=a\mod b$

\begin{lemma}\label{lem-ls3} $\mu $ is invariant under  $\rho_{u}$, for any   $u>0$. 
\end{lemma}

{\bf  Proof of Lemma \ref{lem-ls3} } Let   $0< t_{1}\leq \cdots \leq t_{k-1}\leq t_{k} < \ff$ and let $f_{1},\ldots, f_{k}$ be bounded  Borel measurable functions  on $S\cup \De$ with $f_{j}(\De)=0$, $j=1,\ldots,k$. Fix some $t$ and $u$.

Since $f_{j}(\De)=0$, $j=1,\ldots,k$, 
\begin{equation}
\prod_{j=1}^{k}f_{j}(X_{t_{j}}\circ \rho_{u} \circ k_{t})=1_{\{t_{k}<t\}}\prod_{j=1}^{k}f_{j}(X_{(t_{j}+u)_{t}})
.\label{10.5x}
\end{equation}
Set $s_{j}=t_{j}+u$. Since $0< t_{1}\leq \cdots \leq t_{k-1}\leq t_{k} < t$, it is clear that for some $i$  and some $l$
\begin{equation}
0\leq s_{l} -it\leq \cdots \leq s_{k} -it \leq s_{1}-(i-1)t \ldots \leq s_{l-1}-(i-1)t< t\label{10.5xz}
\end{equation}
Therefore, by (\ref{br7})  
\bea
&&
Q_{t}^{x,x}\(   \prod_{j=1}^{k}f_{j}(X_{t_{j}})\circ \rho_{u}\circ k_{t}\) \label{ls.4q}\\ 
&& =  1_{\{t_{k}<t\}}\int_{S_{k}} p_{s_{l}-it}(x,y_{l}) f_{l}(y_{l}) p_{t_{l+1}-t_{l}}(y_{l},y_{l+1})f_{l+1}(y_{2})\cdots \nonumber\\
&&\hspace{.3 in}\cdots   p_{t_{k} -t_{k-1} }(y_{k-1},y_{k})f_{k}(y_{k}) p_{t_{1}+t-t_{k}}(y_{k},y_{1})f_{1}(y_{1})\cdots \nonumber\\
&&\hspace{.3 in}\cdots  p_{t_{l-1}-t_{l-2}}(y_{l-2},y_{l-1})f_{l-1}(y_{l-1})p_{it-s_{l-1}}(y_{l-1},x)\,dm(y_{1})
\cdots \,dm (y_{k}).\nn
\eea
Integrating both sides with respect to $dm (x)$ and applying the Chapman-Kolmogorov equation (\ref{CK})  we obtain
\bea
&&
\int_{S}Q_{t}^{x,x}\(   \prod_{j=1}^{k}f_{j}(X_{t_{j}})\circ \rho_{u}\circ k_{t}\) \,dm (x)\label{ls.4q1}\\ 
&& =  1_{\{t_{k}<t\}}\int_{S_{k}}   f_{l}(y_{l}) p_{t_{l+1}-t_{l}}(y_{l},y_{l+1})f_{l+1}(y_{2})\cdots \nonumber\\
&&\hspace{.3 in}\cdots   p_{t_{k} -t_{k-1} }(y_{k-1},y_{k})f_{k}(y_{k}) p_{t_{1}+t-t_{k}}(y_{k},y_{1})f_{1}(y_{1})\cdots \nonumber\\
&&\hspace{.3 in}\cdots  p_{t_{l-1}-t_{l-2}}(y_{l-2},y_{l-1})f_{l-1}(y_{l-1})p_{t_{l}-t_{l-1}}(y_{l-1},y_{l})\,dm(y_{1})
\cdots \,dm (y_{k}).\nn
\eea
where in the last line we used 
\[it-s_{l-1}+s_{l} -it=s_{l}-s_{l-1}=t_{l}-t_{l-1}.\]
Comparing with (\ref{ls.4}) we obtain our Lemma.\qed

 \section{The loop soup}\label{sec-soup}

 Let $\Om$ be the path space for $X$ described after (\ref{radon}).
For any $\al>0$, let $\mathcal{L}_{\al}$ be  a Poisson point process  on $\Om$ with intensity measure $\al \mu$.   Note that  $\mathcal{L}_{\al}$ is a random variable; each realization  of $\mathcal{L}_{\al}$ is countable subset of $\Om$. To be more specific, let\begin{equation}
N(A):=\# \{\mathcal{L}_{\al}\cap A \},\hspace{.2 in}A\subseteq \Om.\label{po.1}
\end{equation}
Then for any disjoint measurable subsets $A_{1},\ldots ,A_{n}$ of $\Om$, the random variables $N(A_{1}),\ldots ,N(A_{n})$, are independent, and $N(A)$ is a Poisson random variable with parameter $\al \mu (A)$, i.e.
\begin{equation}
P\(N(A)=k\)={(\al \mu (A))^{k} \over k!}e^{-\al \mu (A)}.\label{po.2}
\end{equation}
%added
(When $\mu(A)=\infty$, this means that $P(N(A) =\infty) =1$.)
%end added
  The Poisson point process $\mathcal{L}_{\al}$ is called the `loop soup' of the Markov process $X$. 
 The term `loop soup' is used in \cite{LW}, \cite{LF} and \cite[Chapter 9]{LL}. In \cite{Le Jan} $\mathcal{L}_{\al}$ is referred to, less colorfully albeit more descriptively, as 
 Poissonian ensembles of Markov loops. See also \cite{Sz1} and \cite{Lupu}.
 
 We define
the `loop soup local time', $\wh L^{x}$, of $X$, by
\begin{equation}
\wh L ^{x}=\sum_{\om\in \mathcal{L_{\al}}}L_{\ff}^{x}(\om). \label{ls.10}
\end{equation} 
 
\medskip	 The next theorem is given for associated Gaussian squares in \cite[Theorem 9]{Le Jan1}.  

\begin{theorem}\label{theo-9.1}    Let   $X$ be a transient Borel right process with state space $S$, as described in the beginning of this section, and let $u(x,y)$, $x,y\in S$ denote its  potential density. Let  $\{ \wh L ^{x} ,x\in S \}$   be the   loop soup local time   of $X$. Then   $\{\,\wh L ^{x} ,x\in S \}$, is an  $\al$-permanental process with kernel $u(x,y)$.  
 \end{theorem}
 
  \Proof     
  By the master formula for Poisson processes, \cite[(3.6)]{K},
\begin{equation}
E\(e^{ \sum_{j=1}^{n}z_{j}\wh L ^{x_{j}} }\)=\exp \(\al\(\int_{\Om}\( e^{\sum_{j=1}^{n}z_{j}L_{\ff}^{x_{j}}(\om)}-1\)\,d\mu(\om)  \)\).\label{ls.11}
\end{equation}
Differentiating  each side of (\ref{ls.11}) with respect to $z_{1},\ldots, z_{n}$ and then setting $z_{1},\ldots, z_{n}$ equal to zero, we get
\begin{equation}
E\( \prod_{j=1}^{n}\wh L ^{x_{j}}  \)=\sum_{l=1} ^{n}\,\,\sum_{\cup_{i=1}^{l} B_{i}=[1,n]}\,\,\al^{l}\,\prod_{i=1}^{l}\,\,\mu\(\prod_{j\in B_{i}} L_{\ff}^{x_{j}}\)\label{mom.11},
\end{equation}
where the second  sum is over all partitions $B_{1},\ldots, B_{l}$ of $[1,n]$. Using (\ref{ls.5c4}) it is easily seen that this is 
\be
E\( \prod_{j=1}^{n}\wh L ^{x_{j}}  \)=  \sum_{\pi\in \mathcal{P}_{n}}\al^{c(\pi)} \prod_{j=1}^{n}u (x_{j},x_{\pi(j)})\label{ls.14f}.
\ee
\qed

 Let $ \th=\{\th_{x} , x\in S\}$ be an $\al$-permanental process   with kernel $u(x,y)$, $x,y\in S$, as considered in Theorem \ref{theo-9.1}. Clearly, by our loop soup construction,  $\th$ is infinitely divisible. In  \cite[Corollary 3.4]{EK}, Eisenbaum and Kaspi show that   the
L\'{e}vy measure of $ \{{\th_{x} }, x\in S\}$  is given by the law of  $\{L_{\ff}^{x}, x\in S\}$ under the $\si$-finite measure 
\be
 {\al  \over L_{\ff}^{y}}Q^{y ,y}  
\ee
for any $y\in S$.  
However it follows from Theorem \ref{theo-9.1} that the  loop measure $\al\mu$ is the L\'evy measure of  $ \{{\th_{x} }, x\in S\}$. Therefore
\be
\{L_{\ff}^{x}, x\in S; \mu\}\stackrel{law}{=} \{L_{\ff}^{x}, x\in S; { 1\over L_{\ff}^{y}}Q^{y ,y} \},
 \ee
 for any $y\in S$. This fact is also an immediate consequence of Lemma \ref{lem-qmu}.

 \section{The Isomorphism Theorem via loop soup}\label{sec-LIT}
 
 For our Isomorphism Theorem we will need a special case of the Palm formula for Poisson processes $\mathcal{L}$ with intensity measure $n$ on a measurable space $\mathcal{S}$, see \cite[Lemma 2.3]{Bertoin}. This says that for any positive function $f$ on $\mathcal{S}$ and any measurable functional $G$ of $\mathcal{L}$
 \begin{equation}
E_{\mathcal{L}} \(\(\sum_{\om\in \mathcal{L}}f(\om)   \)G(\mathcal{L})\)=
\int E_{\mathcal{L}} \(G(\om'\cup\,\mathcal{L})\)f(\om')\,dn(\om'). \label{palm.0}
 \end{equation}

\bt[Isomorphism Theorem]  \label{theo-lilt}  For any $x, x_{1},x_{2},\ldots \in S$ and any  bounded 
measurable function  
$F$ on $R^\ff_+$,  
\begin{equation}
E_{\mathcal{L}_{\al}}  Q^{x,x}\(  F \( \wh L^{x_{j}} +L^{x_{j}}\) \) ={1 \over \al}E_{\mathcal{L}_{\al}} \(\wh L^{x} \,F\(\wh L^{x_{j}}\)\).\label{6.1}
\end{equation} 
  (Here we use the notation $F( f(x_{j})):=F( f(x_{ 1}),f(x_{ 2}),\ldots 
)$.)
\et 

\Proof  
We apply the Palm formula with intensity measure $\al\mu$, 
\begin{equation}
f(\om)=  L^{x}(\om)\label{palm.2}
\end{equation}
and 
\begin{equation}
G(\mathcal{L})=F\( \wh L^{x_{j}}\).\label{palm.2a}
\end{equation}
Note that 
\begin{equation}
\sum_{\om\in \mathcal{L}}f(\om)=\wh L^{x}.\label{palm.2aa}
\end{equation}
Also
\bea
&&
 \wh L^{x_{j}}(\om'\cup\mathcal{L_{\al}})
=\sum_{\om\in \mathcal{\om' \cup L_{\al}}}  L^{x_{j}} (\om)\nn\\ 
&&=\wh L^{x_{j}}(\mathcal{L_{\al}})+L^{x_{j}}(\om'),\label{palm.2b}
\eea
so that 
\begin{equation}
G(\om'\cup\mathcal{L_{\al}})=F\(\wh L^{x_{j}}(\mathcal{L_{\al}})+L^{x_{j}}(\om')\).\label{palm.3}
\end{equation}
Then by (\ref{palm.0})
 \be 
E_{\mathcal{L}_{\al}} \(\wh L^{x} \,F\(\wh L^{x_{j}}\)\) \label{palm.4} 
=\al E_{\mathcal{L}_{\al}}  \mu\(L^{x} \,  F \( \wh L^{x_{j}} +L^{x_{j}}\) \),
 \ee 
 so that our Theorem follows from Lemma \ref{lem-qmu}.\qed

  \section{Invariance under loop rotation}\label{sec-inv}
  
\def\ov{\overline}
\def\mod{{\mathop{\rm mod}}}

In subsection \ref{sec-rot}, assuming the existence of transition densities, we  gave a simple proof of the fact that the loop measure is invariant under loop rotation.   In this section we give a proof of this fact without assuming transition densities. This proof is considerably more complicated.

Because the lifetime $\zeta$ is rotation invariant ($\zeta(\rho_v\omega) =\zeta(\omega)$ so long as $\zeta(\omega)<\infty$), the rotation invariance of the loop measure $\mu$ is equivalent to that of the measure $\nu$ defined by $\nu(F) :=\mu(\zeta F)$. By (\ref{6.41}) and (\ref{ls.3n}) we have 
\begin{equation}
\nu(F)=\int_{S}Q^{x,x}\( F \)\,dm (x)= \int_{S} P^{x}\(\int_{0}^{\ff} \,F\circ k_{t}\, \,d_{t}L_{t}^{x}\)\,dm (x).\label{nuq1}
\end{equation}

The measure $\nu$ is more convenient for the calculations that follow, because of the following formula, where $\Gamma_{x,y}$ is defined in (\ref{gamdef}):

\bl\label{lem-fi} For a.e. $0<t_1<t_2<\cdots<t_k$,
\bea
&&
\nu(X_{t_1}\in dy_1,\ldots,X_{t_k}\in dy_k, \zeta\in dt)\label{p00}\\
&&\hspace{.5 in} = m(dy_1)\prod_{j=2}^k P_{t_j-t_{j-1}}(y_{j-1},dy_j)\,\Gamma_{y_k,y_1}(dt-t_k+t_1)
\nn
\eea
 as measures on the product space $S^k\times(t_k,\infty)$. Furthermore, with $t_{1}=0$, (\ref{p00}) holds for all $0<t_2<t_3<\cdots<t_k$.
\el

{\bf  Proof of Lemma \ref{lem-fi}: }  Using  (\ref{nuq1}) we see that
\be
\nu\(\prod_{j=1}^{k}f_{j}(X_{t_{j}})e^{-\bb\ze}\)  
\label{nuq2} =\int_{S}P^{x}\(  \prod_{j=1}^{k}f_{j}(X_{t_{j}})    \int_{t_{k}}^{\ff} \,e^{-\bb t} \,d_{t}L_{t}^{x}\)\,dm (x).    
\ee
Using the Markov property and (\ref{res.3}) we see that
\begin{eqnarray}
&&P^{x}\(  \prod_{j=1}^{k}f_{j}(X_{t_{j}})    \int_{t_{k}}^{\ff} \,e^{-\bb t} \,d_{t}L_{t}^{x}\)
\label{nuq3}\\
&&= P^{x}\(  \prod_{j=1}^{k}f_{j}(X_{t_{j}}) e^{-\bb t_{k}}   \( \int_{0}^{\ff} \,e^{-\bb t} \,d_{t}L_{t}^{x}\)  \circ \th _{t_{k}}    \)  \nonumber\\
&&= P^{x}\(  \prod_{j=1}^{k}f_{j}(X_{t_{j}}) e^{-\bb t_{k}}  E^{X_{t_{k}}} \( \int_{0}^{\ff} \,e^{-\bb t} \,d_{t}L_{t}^{x}\)     \)  \nonumber\\
&&= P^{x}\(  \prod_{j=1}^{k}f_{j}(X_{t_{j}}) e^{-\bb t_{k}} \,\, u^{\bb}\(X_{t_{k}},x\)\) \nonumber \\
&&=\int_{S^{k}} P^{\bb}_{t_1}(x,dy_1)f_{1}(y_1)\,\(\prod_{j=2}^k P^{\bb}_{t_j-t_{j-1}}(y_{j-1},dy_j)f_{j}(y_j)\)\,\, u^{\bb}\(y_{k},x\). \nonumber
\end{eqnarray}
Here $P^{\bb}_{t }(x,dy )=e^{-\bb t}P_ {t }(x,dy )$. 
Using (\ref{intp}) and then the Markov property as in the previous display
\begin{eqnarray}
&& \int_{t_{1}}^{\ff} \,e^{-\bb t}\Ga_{y_k,y_1} (dt)
\label{res.4}\\
&&=E^{y_k}\(  \int_{t_{1}}^{\ff} \,e^{-\bb t} \,d_{t}L_{t}^{y_{1}} \)   \nonumber\\
&&=e^{-\bb t_{1}}\,E^{y_k}\(  \( \int_{0}^{\ff} \,e^{-\bb t} \,d_{t}L_{t}^{y_{1}}\)  \circ \th _{t_{1}} \)\nonumber\\
&& =e^{-\bb t_{1}} \,E^{y_k}\(u^{\bb}\(X_{t_{1}}, y_{1}\)\)\nonumber\\
&& =\int_{S} P^{\bb}_{t_1}(  y_k,\,dz)\,u^{\bb}\(z, y_{1}\).\nonumber
\end{eqnarray}
We claim that for a.e. $t_{1}$, as measures in $y_{1}$,
\be
\int_{x\in S}  u^{\bb}\(y_{k},x\) P^{\bb}_{t_1}(x,dy_1)    \,dm (x)=\int_{z\in S} P^{\bb}_{t_1}(  y_k,\,dz)\,u^{\bb}\(z, y_{1}\)\,dm (y_{1}).\label{res.5}
\ee 
To see this, it suffices to integrate both sides with respect to $e^{-\al t_{1}}\,dt_{1}$, use (\ref{res.2}) with $\al$ replaced by  $\al+\bb$, and the fact that $S$ has a countable base. (It is important to note that we allow $y_{k}=y_{1}$).

Combining (\ref{nuq2})-(\ref{res.5}) we obtain for a.e. $t_{1}$
\bea
&&
\nu\(\prod_{j=1}^{k}f_{j}(X_{t_{j}})e^{-\bb\ze}\)\label{nuq5}\\
&& =\int f_{1}(y_1)\,\prod_{j=2}^k P^{\bb}_{t_j-t_{j-1}}(y_{j-1},dy_j)f_{j}(y_j)\,\,\int_{t_{1}}^{\ff} \,e^{-\bb t}\Ga_{y_k,y_1} (dt) \, \,dm (y_{1}).
\nn
\eea
This agrees with what we obtain from the right hand side of (\ref{p00}):
\begin{eqnarray}
&&\int_{S^k\times(t_k,\infty)}     m(dy_1)\prod_{j=2}^k P_{t_j-t_{j-1}}(y_{j-1},dy_j)\,\Gamma_{y_k,y_1}(dt-t_k+t_1)  
\label{nuk6}\\
&&   \hspace{3 in}\(\prod_{j=1}^{k}f_{j}(y_{j})e^{-\bb t}\) \nonumber\\
&&=\int_{S^k}     m(dy_1)f_{1}(y_{1})\prod_{j=2}^k P_{t_j-t_{j-1}}(y_{j-1},dy_j)f_{j}(y_{j})\, 
\nn\\
&&   \hspace{2 in}\int_{t_{k}}^{\ff}e^{-\bb t}\Gamma_{y_k,y_1}(dt-t_k+t_1)  \nonumber\\
&&=\int_{S^k}     m(dy_1)f_{1}(y_{1})\prod_{j=2}^k P_{t_j-t_{j-1}}(y_{j-1},dy_j)f_{j}(y_{j})\, 
\nn\\
&&   \hspace{2 in}e^{-\bb (t_k-t_1)}\int_{t_{1}}^{\ff}e^{-\bb t}\Gamma_{y_k,y_1}(dt).  \nonumber
\end{eqnarray}
This completes the proof of our Lemma when $t_{1}>0$.

When $t_{1}=0$, it follows from (\ref{nuq2})-(\ref{nuq3}), and then (\ref{res.3}) and (\ref{intp}) that 
\begin{eqnarray}
&&
\nu\(\prod_{j=1}^{k}f_{j}(X_{t_{j}})e^{-\bb\ze}\)  
\label{nuq27}\\
&&=\int_{S} f_{1}(x)\,\prod_{j=2}^k P^{\bb}_{t_j-t_{j-1}}(y_{j-1},dy_j)f_{j}(y_j)\,\, u^{\bb}\(y_{k},x\) \,dm (x)  \nonumber \\
&&=\int_{S} f_{1}(x)\,\prod_{j=2}^k P^{\bb}_{t_j-t_{j-1}}(y_{j-1},dy_j)f_{j}(y_j)\,\,\int_{0}^{\ff} \,e^{-\bb t}\Ga_{y_k,x} (dt) \,\,dm (x).   \nonumber
\end{eqnarray}
This agrees with (\ref{nuq5}) for $t_{1}=0$, and the rest of the proof follows as in (\ref{nuk6}).
\qed

As a byproduct of our proof we now show that
\begin{equation}
\sup_{t_{1}\geq 0}\,\nu\(f_{1}(X_{t_{1}})e^{-\bb\ze}\)<\ff,\label{dom.1}
\end{equation}
 for any continuous compactly supported $f_{1}$.
To see this, note that by (\ref{nuq2})-(\ref{nuq3})
\begin{equation}
\nu\(f_{1}(X_{t_{1}})e^{-\bb\ze}\)=\int_{S} \int_{S} P^{\bb}_{t_1}(x,dy_1)f_{1}(y_1)\,  u^{\bb}\(y_{1},x\)  \,dm (x).  \label{dom.2}
\end{equation}
By (\ref{res.5}), for a.e. $t_{1}$ this equals
\begin{equation}
\int_{S}\(\int_{S} P^{\bb}_{t_1}(  y_1,\,dz)\,u^{\bb}\(z, y_{1}\)\)f_{1}(y_1)\,dm (y_{1}).\label{dom.3}
\end{equation}
But as noted in the paragraph containing (\ref{res.1}),  $u^{\bb}\(z, y_{1}\)$ is bounded, uniformly in $z$ for $y_{1}$ in the compact support of $f_{1}(y_1)$. Hence (\ref{dom.3}) is bounded by
\begin{equation}
C\int_{S}\(\int_{S} P^{\bb}_{t_1}(  y_1,\,dz)\) f_{1}(y_1)\,dm (y_{1})\leq C\int_{S}f_{1}(y_1)\,dm (y_{1}).  \label{dom.4}
\end{equation}
Thus we have shown that for some dense $D\subseteq R_{+}^{1}$
\begin{equation}
\sup_{t_{1}\in D}\nu\(f_{1}(X_{t_{1}})e^{-\bb\ze}\)\leq C\int_{S}f_{1}(y_1)\,dm (y_{1}),\label{dom.1v}
\end{equation}
and the right hand side is finite by (\ref{radon}).
 (\ref{dom.1}) then follows using right continuity.

We will also need the following.
\bl\label{lem-pg}
\begin{equation}
P_s(x,dy)\,ds=m(dy)\Gamma_{x,y}(ds).\label{pg1}
\end{equation}
\el

{\bf  Proof of Lemma \ref{lem-pg}: }We have 
\begin{eqnarray}
&&\int \Gamma_{x,y}([0,t])\,f(y)\,m(dy)
\label{pg2}\\
&& =\int E^{x}\(L_{t}^{y}\)\,f(y)\,m(dy) =\int E^{x}\(L_{\ff}^{y}-L_{\ff}^{y}\circ \th_{t}\)\,f(y)\,m(dy) \nonumber\\
&& = \int \(u(x,y)-E^{x}\(u(X_{t},y)\)\)\,f(y)\,m(dy) \nonumber\\
&& =\int \int_{0}^{t} P_s(x,dy)\,f(y)\,ds. \nonumber
\end{eqnarray}
\qed

Let us define the process $\overline X$ to be the periodic extension of $X$; that is,
\begin{equation}
\ov X_t
	=\left\{\begin{array}{ll}
X_{t-q\zeta}, & \mbox{ if $q\zeta\le t<(q+1)\zeta$, $q=0,1,2,\ldots$}\\
	 X_t, & \mbox{if $\zeta=\infty$}
	\end{array} \right.
	\label{p2}
\end{equation} 

It will be convenient to write
\be
\ov I_\alpha(f):=\int_0^\infty e^{-\alpha t} f(\ov X_t)\,dt,\qquad   I_\alpha(f):=\int_0^\infty e^{-\alpha t} f(X_t)\,dt.
\label{p3}
\ee
The key observation is that
\be
\ov I_\alpha(f) ={I_\alpha(f)\over 1-e^{-\alpha \zeta}},
\label{p4}
\ee
for all $\alpha>0$. This follows from
\begin{eqnarray}
\ov I_\alpha(f):&=&\int_0^\infty e^{-\alpha t} f(\ov X_t)\,dt
\nn\\
&=& \sum_{q=0}^{\ff} \int_{q\ze}^{(q+1)\ze} e^{-\alpha t} f(\ov X_t)\,dt \nonumber\\
&=& \sum_{q=0}^{\ff}e^{-\alpha q\ze}  \int_{0}^{\ze} e^{-\alpha t} f( X_t)\,dt ={I_\alpha(f)\over 1-e^{-\alpha \zeta}}\,.\nonumber
\end{eqnarray}
Hence for any continuous compactly supported $f$
\bea
&&
\int_0^\infty e^{-\alpha t}\nu\(\(1-e^{-\alpha \zeta}\)f(\ov X_t)e^{-\bb\ze}\) \,dt\label{subst.1}\\
&&=\nu\(\(1-e^{-\alpha \zeta}\) \ov I_\alpha(f)e^{-\bb\ze}\)
=\nu\(I_\alpha(f)e^{-\bb\ze}\)\nn\\
&&=\int_0^\infty e^{-\alpha t}\nu\( f( X_t)e^{-\bb\ze}\) \,dt<\ff\nn
\eea
by (\ref{dom.1}). It follows that for any $\al $
\begin{equation}
\nu\(\(1-e^{-\alpha \zeta}\)f(\ov X_t)e^{-\bb\ze}\)<\ff, \mbox{ for a.e. $t$.}\label{dom.10}
\end{equation}

The rotation invariance of $\mu$ or $\nu$ is equivalent to the statement that
\be
\nu\left(\prod_{j=1}^k f_j(\ov X_{t_j+r})1_{\{t_{k}<\ze\}}\right) =
\nu\left(\prod_{j=1}^k f_j(\ov X_{t_j})1_{\{t_{k}<\ze\}}\right)
 \label{p18glm}
\ee 
for all  $0<t_1<\cdots< t_k$ and $r>0$ and all $f_{j}\geq 0$ continuous with compact support.
This will follow once we show that the joint distribution of $(\ov X, \zeta)$ is invariant under time shifts.  That is, $((\ov X_{t+v})_{t\ge 0}, \zeta)$ has the same distribution (under $\nu$)  as $((\ov X_t)_{t\ge 0},\zeta)$ for all $v>0$.

To prove this we will first show that for all $k$ and all $\al_{1}, \ldots, \al_{k},$
\bea
&&
\int_{[0,\ff)^{k}} e^{-\sum_{j=1}^{k}\al_{j}t_{j}}\,\,\nu\left(\prod_{j=1}^k f_j(\ov X_{t_j})\,g(\ze)\right)\prod_{j=1}^k\,dt_{j} \label{p18}\\
&&\hspace{.4 in} =\int_{[0,\ff)^{k}} e^{-\sum_{j=1}^{k}\al_{j}t_{j}}\,\,F_{k}(t_1,\ldots,t_k)\prod_{j=1}^k\,dt_{j}
\nn
\eea
for all $g$ of the form $g(\ze)=(1-e^{-\al \zeta})e^{-\beta \zeta}$, and where
\bea
&&
F_{k}(t_1,\ldots,t_k) = \sum_{\sigma\in\PP_k}1_{\{0\leq t_{\si (1)}\leq\cdots\leq t_{\si (k)}\}}\label{p18t}\\
&&
\hspace{1.5 in}\nu\left(f_{\si (1)}(X_0)\prod_{j=2}^k f_{\si (j)}(\ov X_{t_{\si (j)}-t_{\si (1)}})\,g(\ze)\right).
\nn
\eea
By (\ref{subst.1}) the left hand side of (\ref{p18}) is finite for all $\al_{1}\geq \al$, while the right hand side is finite since
\begin{equation}
\nu\(f_{j}(X_0)\)=\int_{S}Q^{x,x}\( f_{j}(X_0) \)\,dm (x)=\int_{S}u(x,x)f_{j}(x)\,dm (x)<\ff.\label{fino}
\end{equation}
By uniqueness of Laplace transforms, it 
then
follows that
\bea
&&
\nu\left(\prod_{j=1}^k f_j(\ov X_{t_j})\,g(\ze)\right) \label{p18g}
 =F_{k}(t_1,\ldots,t_k)
\eea
for Lebesgue a.e.\ $k$-tuple $(t_1,\ldots,t_k)$, and in particular, for any $r>0$,
\bea
&&
\int_{[r,\ff)^{k}} e^{-\sum_{j=1}^{k}\al_{j}t_{j}}\,\,\nu\left(\prod_{j=1}^k f_j(\ov X_{t_j})\,g(\ze)\right)\prod_{j=1}^k\,dt_{j} \label{p18m}\\
&&\hspace{.4 in} =\int_{[r,\ff)^{k}}  e^{-\sum_{j=1}^{k}\al_{j}t_{j}}\,\,F_{k}(t_1,\ldots,t_k)\prod_{j=1}^k\,dt_{j}.
\nn
\eea
It follows that for any $r>0$,
\bea
&&
\int_{[0,\ff)^{k}} e^{-\sum_{j=1}^{k}\al_{j}(t_{j}+r)}\,\,\nu\left(\prod_{j=1}^k f_j(\ov X_{t_j+r})\,g(\ze)\right)\prod_{j=1}^k\,dt_{j} \label{p18k}\\
&&\hspace{.4 in} =\int_{[0,\ff)^{k}} e^{-\sum_{j=1}^{k}\al_{j}(t_{j}+r)}\,\,F_{k}(t_1+r,\ldots,t_k+r)\prod_{j=1}^k\,dt_{j}.
\nn
\eea
But it is easily seen that $F_{k}(t_1+r,\ldots,t_k+r)=F_{k}(t_1,\ldots,t_k)$ so that, canceling the common constant factor $e^{-\sum_{j=1}^{k}\al_{j}r}$, we obtain
\bea
&&
\int_{[0,\ff)^{k}} e^{-\sum_{j=1}^{k}\al_{j}t_{j}}\,\,\nu\left(\prod_{j=1}^k f_j(\ov X_{t_j+r})\,g(\ze)\right)\prod_{j=1}^k\,dt_{j} \label{p18n}\\
&&\hspace{.4 in} =\int_{[0,\ff)^{k}} e^{-\sum_{j=1}^{k}\al_{j}t_{j}}\,\,F_{k}(t_1,\ldots,t_k)\prod_{j=1}^k\,dt_{j},
\nn
\eea
and thus comparing with (\ref{p18}) we have that for each $r>0$
\bea
&&
\int_{[0,\ff)^{k}} e^{-\sum_{j=1}^{k}\al_{j}t_{j}}\,\,\nu\left(\prod_{j=1}^k f_j(\ov X_{t_j+r})\,g(\ze)\right)\prod_{j=1}^k\,dt_{j} \label{p18p}\\
&&\hspace{.4 in} =\int_{[0,\ff)^{k}} e^{-\sum_{j=1}^{k}\al_{j}t_{j}}\,\,\nu\left(\prod_{j=1}^k f_j(\ov X_{t_j})\,g(\ze)\right)\prod_{j=1}^k\,dt_{j} .
\nn
\eea
It follows that 
\begin{equation}
\nu\left(\prod_{j=1}^k f_j(\ov X_{t_j+r})\,g(\ze)\right)=\nu\left(\prod_{j=1}^k f_j(\ov X_{t_j})\,g(\ze)\right), \hspace{.1 in}\mbox{  for a.e.  $t_1,\ldots,t_k$.}\label{res.7}
\end{equation} 
This holds for any $k$, in particular for $k=1$, so that using (\ref{dom.10}) we have
\begin{equation}
\nu\left( f_1(\ov X_{t_1+r})\,g(\ze)\right)=\nu\left( f_1(\ov X_{t_1})\,g(\ze)\right)<\ff, \hspace{.2 in}\mbox{  for a.e.  $t_1$.}\label{res.7a}
\end{equation} 
Thus by Fubini we can find a set $T\subseteq R_{+}$ with $T^{c}$ of Lebesgue measure $0$ such that for all $t_{1}\in T$ we have (\ref{res.7a}), and (\ref{res.7}) holds for a.e.  $t_2,\ldots,t_k$. Using the boundedness and  continuity of the $f_{j}$ and the right continuity of $\bar X_{t}$ it follows from the Dominated Convergence Theorem that  (\ref{res.7}) holds for all  $(t_{1},t_2,\ldots,t_k)\in T\times R^{k-1}_{+}$. Let now $f_{1,n}$ be a sequence of continuous functions with compact support with the property that $f_{1,n}\uparrow 1$. By the above, (\ref{res.7}) with $f_{1}$ replaced by $f_{1,n}$ holds for all  $(t_{1},t_2,\ldots,t_k)\in T_{n}\times R^{k-1}_{+}$ for an appropriate $T_{n}\subseteq R_{+}$ with $T_{n}^{c}$ of Lebesgue measure $0$. In particular $T_{\ast}=\cap_{n}T_{n}\not=\emptyset$, and we can apply the Monotone Convergence Theorem with $t_{1}\in T_{\ast}$ to conclude, spelling out $g(\ze)$, that 
\begin{equation}
\nu\left((1-e^{-\al \zeta})\prod_{j=2}^k f_j(\ov X_{t_j+r})e^{-\beta \zeta}\right)=\nu\left((1-e^{-\al \zeta})\prod_{j=2}^k f_j(\ov X_{t_j})e^{-\beta \zeta}\right)\label{res.17}
\end{equation} 
for all $t_2,\ldots,t_k$. Applying once again the Monotone Convergence Theorem for $\al\to\ff$ we obtain
\begin{equation}
\nu\left(\prod_{j=2}^k f_j(\ov X_{t_j+r})e^{-\beta \zeta}\right)=\nu\left(\prod_{j=2}^k f_j(\ov X_{t_j})e^{-\beta \zeta}\right)\label{res.17d}
\end{equation} 
for all $t_2,\ldots,t_k$. Fix a compact $K\subseteq S$. If we replace $f_{2}$ by a sequence $f_{2,n}\uparrow 1_{K}$ and then set   $t_{2}=0$,  we can conclude from (\ref{res.17d}) and (\ref{fino}) that  the finite measures $1_{K}(X_{0})\cdot \nu$ and $1_{K}(X_{0})\cdot\rho_{r\,\ast} \nu$ agree on the $\si$-algebra generated by $\bar X_{t}, t\geq 0$ and $\ze$. Since this holds for any compact $K\subseteq S$, so do 
$\nu$ and $\rho_{r\,\ast} \nu$. Here and below we use the notation  $f_{\ast}\nu (A)=\nu(f^{-1}(A))$.

It remains to prove (\ref{p18}).
Using (\ref{p4})
\bea
&&
\int e^{-\sum_{j=1}^{k}\al_{j}t_{j}}\,\,\nu\left(\prod_{j=1}^k f_j(\ov X_{t_j})\,g(\ze)\right)\prod_{j=1}^k\,dt_{j}\label{p5}\\
&&=\nu\left(\prod_{j=1}^k \ov I_{\alpha_{j}}(f_{j})\,g(\ze)\right)\nonumber\\
&&= \nu\left(\int e^{-\sum_{j=1}^{k}\al_{j}t_{j}}\prod_{j=1}^k {f_j( X_{t_j}) \over 1-e^{-\alpha_{j} \zeta}}\prod_{j=1}^k\,dt_{j}\,\,g(\ze)\right)=\sum_{\sigma\in\PP_k}J_{k}(\si),\nn
\eea
where 
\be
J_{k}(\si):=\nu\left(\int_{0<t_{1}<\cdots<t_{k}} e^{-\sum_{j=1}^{k}\al_{\sigma(j)}t_{j}}\prod_{j=1}^k {f_{\sigma(j)}( X_{t_j}) \over 1-e^{-\alpha_{\sigma(j)} \zeta}}\prod_{j=1}^k\,dt_{j}\,\,g(\ze)\right),
\nn
\ee
and using (\ref{p00})
\bea
&&
J_{k}(\si)=\int_{0<t_{1}<\cdots<t_{k}<t}g(t)\int_{S^k}  
{e^{-\alpha_{\sigma(1)}t_{1}} \over 1-e^{-\alpha_{\sigma(1)} t}}
m(dy_{1}) f_{\sigma(1)}(y_{1})\label{p6}\\
&&\hspace{1 in}  \prod_{j=2}^k {e^{-\alpha_{\sigma(j)}t_{j}} \over 1-e^{-\alpha_{\sigma(j)} t}}
P_{t_{j}-t_{j-1 }}(y_{j-1},dy_{j})f_{\sigma(j)}(y_{j})\,\nn\\
&&\hspace{1.8 in} \Gamma_{y_{k},y_{1}}(dt-t_{k}+t_{1})
\,dt_{1}\cdots dt_{k}.\nn
\eea
We now make the change of variables $r=t_{1}$, $s_{j}=t_{j}-t_{j-1}$ ($j=2,\ldots,k$), $s_{1}=t-t_{k}+t_{1}$ (with accompanying limits of integration $0<r<s_{1}$, $s_{j}>0$) and then integrate out $r$. Writing $\hat s_{j}:= s_{2}+\cdots+ s_{j}$ and $\bar s:= \sum_{j=1}^{k} s_{j}$, the expression in (\ref{p6}) is thereby transformed to
\bea
&&
J_{k}(\si)=\int_{s_{1}>0,\ldots,s_{k}>0} g(\bar s)\int_{S^k}  
\(\int_{0}^{s_{1}}e^{-(\alpha_1+\cdots\alpha_k)r}\,dr\)
m(dy_{1}) f_{\sigma(1)}(y_{1})\nn\\
&&\hspace{.1 in} {1 \over 1-e^{-\alpha_{\sigma(1)}\bar s}}\prod_{j=2}^k {e^{-\alpha_{\sigma(j)} \hat  s_{j}} \over 1-e^{-\alpha_{\sigma(j)}\bar s}}
P_{s_{j}}(y_{j-1},dy_{j})f_{\sigma(j)}(y_{j})\, \Gamma_{y_{k},y_{1}}(ds_{1})
\,ds_{2}\cdots ds_{k}\nn\\
&&=
\(\sum_{j=1}^{k}\alpha_j\)^{-1}
\int_{s_{1}>0,\ldots,s_{k}>0}g(\bar s)\int_{S^k} 
m(dy_{1}) f_{\sigma(1)}(y_{1}){(1-e^{-(\alpha_1+\cdots\alpha_k)s_{1}}) \over 1-e^{-\alpha_{\sigma(1)}\bar s}}\nn\\
&&\hspace{.2 in} \prod_{j=2}^k {e^{-\alpha_{\sigma(j)} \hat  s_{j}} \over 1-e^{-\alpha_{\sigma(j)}\bar s}}
P_{s_{j}}(y_{j-1},dy_{j})\,f_{\sigma(j)}(y_{j})  \Gamma_{y_{k},y_{1}}(ds_{1})
\,ds_{2}\cdots ds_{k}.\label{p7}
\eea 
Using (\ref{pg1}), we can write (\ref{p7})  as
\bea
&&
J_{k}(\si)=\(\sum_{j=1}^{k}\alpha_j\)^{-1}
\int_{s_{1}>0,\ldots,s_{k}>0}g(\bar s)\,
 (1-e^{-(\alpha_1+\cdots\alpha_k)s_{1}})\nn\\
&&\hspace{.4 in} \prod_{j=1}^k {e^{-\alpha_{\sigma(j)}\hat s_{j}} \over 1-e^{-\alpha_{\sigma(j)}\bar s}}\int_{S^k} \prod_{j=1}^k
m(dy_{j})f_{\sigma(j)}(y_{j})
 \Gamma_{y_{j-1},y_{j}}(ds_{j}),
\label{p8}
\eea
where $y_{0}=y_{k}$ and $\hat s_{1}:= 0$.

We now turn to the right hand side of (\ref{p18}) and try to rewrite it in terms which are similar to our last expression for the 
$J_{k}(\si)$'s. Using $\sum_{j=1}^{k}\al_{j}t_{j}=\sum_{j=1}^{k}\al_{\si (j)}t_{\si (j)}$
we have
\begin{eqnarray}
&&R(\al_{1},\ldots,\al_{k}):= \int e^{-\sum_{j=1}^{k}\al_{j}t_{j}}\,\,F_{k}(t_1,\ldots,t_k)\prod_{j=1}^k\,dt_{j}
\label{pt1}\\
&& =\sum_{\sigma\in\PP_k}\int_{\{0\leq t_{\si (1)}\leq\cdots\leq t_{\si (k)}\}}e^{-\sum_{j=1}^{k}\al_{\si (j)}t_{\si (j)}}\nn\\
&&
\hspace{1  in}\nu\left(f_{\si (1)}(X_0)\prod_{j=2}^k f_{\si (j)}(\ov X_{t_{\si (j)}-t_{\si (1)}})\,g(\ze)\right)\prod_{j=1}^k\,dt_{j}  \nonumber\\
&& =\sum_{\sigma\in\PP_k}\int_{\{0\leq t_{1}\leq\cdots\leq t_{k}\}}e^{-\sum_{j=1}^{k}\al_{\si (j)}t_{j}}\nn\\
&&
\hspace{1  in}\nu\left(f_{\si (1)}(X_0)\prod_{j=2}^k f_{\si (j)}(\ov X_{t_{j}-t_{1}})\,g(\ze)\right)\prod_{j=1}^k\,dt_{j}.  \nonumber
\end{eqnarray}
Let us now fix $\sigma\in\PP_k$ and consider  the corresponding term in (\ref{pt1})
\be
\int_{0<t_{1}<\cdots<t_{k}}\prod_{j=1}^k e^{-\alpha_{\sigma(j)}t_{j}}dt_{j}\,
\nu\left(f_{\sigma(1)}(\ov X_0)\prod_{j=2}^kf_{\sigma(j)}(\ov X_{t_{j}-t_{1}})\,g(\ze)\right).
\label{p9}
\ee
Changing variables ($r_{1}=t_{1}$, $r_{j}=t_{j}-t_{1}$ for $j=2,\ldots,k$) and integrating out $r_{1}$,   this can be rewritten as
\be
\(\sum_{j=1}^{k}\alpha_j\)^{-1}
\int_{0<r_{2}<\cdots<r_{k}}\prod_{j=2}^k e^{-\alpha_{\sigma(j)}r_{j}} dr_{j}\,
\nu\left(f_{\sigma(1)}(\ov X_0)\prod_{j=2}^k f_{\sigma(j)}(\ov X_{r_{j}})\,g(\ze)\right).
\label{p10}
\ee

Summing first over all  permutations $\sigma\in\PP_k$ with $\sigma(1) =i$ and then over $i$
we obtain 
\be
R(\al_{1},\ldots,\al_{k})=\(\sum_{j=1}^{k}\alpha_j\)^{-1}\sum_{i=1}^{k}\,\,
\nu\left(f_{i}(\ov X_0)\prod_{j\not=i} \ov I_{\alpha_j}(f_{j})\cdot \,g(\ze)\right).
\label{p11}
\ee
Using (\ref{p4}) we can express this as
\be
R(\al_{1},\ldots,\al_{k})=\(\sum_{j=1}^{k}\alpha_j\)^{-1}\sum_{i=1}^{k}\,\,
\nu\left(f_{i}(X_0)\prod_{j\not=i} { I_{\alpha_{j}}(f_{j})\over 1-e^{-\alpha_{j}\zeta}}\cdot \,g(\ze)\right).
\label{p12}
\ee

Using Lemma \ref{lem-fi} we then see that 
\begin{equation}
R(\al_{1},\ldots,\al_{k})=\sum_{\sigma'\in\PP_k}K_{k}(\si')\label{p112}
\end{equation}
where
\bea
&&
K_{k}(\si'):=\(\sum_{j=1}^{k}\alpha_j\)^{-1}\int_{0<t_{2}<\cdots<t_{k}<t}g(t)\,
\int_{S^k} m(dy_1)f_{\sigma'(1)}(y_1)\nn\\
&&\hspace{1.3 in}  \prod_{j=2}^k P_{t_{j}-t_{j-1}}(y_{j-1},dy_{j})
f_{\sigma'(j)}(y_{j}) \nn\\
&& \hspace{1.4 in} {e^{-\alpha_{\sigma'(j)}t_{j}}\over 1-e^{-\alpha_{\sigma'(j)}t}}\,\Gamma_{y_{k},y_{1}}(dt-t_{k})\,dt_{2}\cdots dt_{k},\label{p13}
\eea
with the convention that $t_{1}=0$. Once more making the change of variables $s_{1}=t-t_{k}$, $s_{2}=t_{2}$, $s_{j}=t_{j}-t_{j-1}$ for $j=3,\ldots,k$, (\ref{p13}) becomes
\bea
&&
K_{k}(\si')=\(\sum_{j=1}^{k}\alpha_j\)^{-1}\int_{s_{1}>0,\ldots,s_{k}>0}\,g(\bar s)
\int_{S^k} m(dy_1)f_{\sigma'(1)}(y_1)\label{p13r}\\
&&\hspace{.5 in}  \prod_{j=2}^k P_{s_{j}}(y_{j-1},dy_{j})
f_{\sigma'(j)}(y_{j}) {e^{-\alpha_{\sigma'(j)}\hat s_{j}}\over 1-e^{-\alpha_{\sigma'(j)}\bar s}}\,\Gamma_{y_{k},y_{1}}(ds_{1})\,ds_{2}\cdots ds_{k}.\nn
\eea
Using (\ref{pg1}) again, we can write (\ref{p13r})  as 
\bea
&&
K_{k}(\si')=\(\sum_{j=1}^{k}\alpha_j\)^{-1}\int_{s_{1}>0,\ldots,s_{k}>0}\,g(\bar s)
\int_{S^k} m(dy_1)f_{\sigma'(1)}(y_1)\label{p13s}\\
&&\hspace{.8 in} \( \prod_{j=2}^k   m(dy_{j})\Gamma_{y_{j-1},y_{j}}(ds_{j})   {e^{-\alpha_{\sigma'(j)}\hat s_{j}}\over 1-e^{-\alpha_{\sigma'(j)}\bar s}}\)
f_{\sigma'(j)}(y_{j})  \nn\\
&& \hspace{2.5 in}  \,\Gamma_{y_{k},y_{1}}(ds_{1})\,ds_{2}\cdots ds_{k}\nn\\
&&
=\(\sum_{j=1}^{k}\alpha_j\)^{-1}\int_{s_{1}>0,\ldots,s_{k}>0}\,g(\bar s)
\int_{S^k} 
(1-e^{-\alpha_{\sigma'(1)}\bar s})\prod_{j=1}^k m(dy_{j})\nn\\
&&\hspace{1.2 in} \Gamma_{y_{j-1},y_{j}}(ds_{j})  f_{\sigma'(j)}(y_{j}) {e^{-\alpha_{\sigma'(j)}\hat s_{j}}\over 1-e^{-\alpha_{\sigma'(j)}\bar s}}, 
 \label{p14}
\eea
where $y_{0}=y_{k}$, and $\hat s_{1} =0$. We reorganize this as
\begin{eqnarray}
&&K_{k}(\si')=\(\sum_{j=1}^{k}\alpha_j\)^{-1}\int_{s_{1}>0,\ldots,s_{k}>0}g(\bar s)
 (1-e^{-\alpha_{\sigma'(1)}\bar s})
\nn\\
&&\hspace{.2 in} \prod_{j=1}^k {e^{-\alpha_{\sigma'(j)}\hat s_{j}} \over 1-e^{-\alpha_{\sigma'(j)}\bar s}}\int_{S^k} \prod_{j=1}^k 
m(dy_{j})f_{\sigma'(j)}(y_{j})
 \Gamma_{y_{j-1},y_{j}}(ds_{j}).  \label{p14reo}
\end{eqnarray}

In view of (\ref{p5}) and (\ref{p112}), to prove  (\ref{p18}) we need to show that    
\begin{equation}
\sum_{\sigma\in\PP_k}J_{k}(\si)=\sum_{\sigma\in\PP_k}K_{k}(\si),\label{pj1}
\end{equation}
and to this end it suffices to show that for each $\sigma^{\ast}\in\PP_k$
\begin{equation}
\sum_{\sigma\sim \sigma^{\ast}}J_{k}(\si)=\sum_{\sigma\sim \sigma^{\ast}}K_{k}(\si),\label{pj2}
\end{equation}
where $\sigma\sim \sigma^{\ast}$ means that $\sigma$ is a `rotation' of $\sigma^{\ast}$. In other words, for some $1\leq l\leq k$ we have \[(\sigma(1),\sigma(2),\ldots, \sigma(k))=(\sigma^{\ast}(l),\sigma^{\ast}(l+1),\ldots, \sigma^{\ast}(k),\sigma^{\ast}(1),\ldots, \sigma^{\ast}(l-1)).\]

 Comparing (\ref{p8}) and (\ref{p14reo}) with $\si'=\si$ we see that the only difference is the presence of $e^{-(\alpha_1+\cdots\alpha_k)s_{1}}$ in (\ref{p8}) while in (\ref{p14reo}),  with $\si'=\si$, this is replaced by $e^{-\alpha_{\sigma(1)}\bar s}$. Thus to prove (\ref{pj2})
  it suffices to show that 
\begin{eqnarray}
&&\int_{s_{1}>0,\ldots,s_{k}>0} g(\bar s)
  e^{-(\alpha_1+\cdots\alpha_k)s_{1}}\label{p8m}\\
&&\hspace{.2 in} \prod_{j=1}^k {e^{-\alpha_{\sigma(j)}\hat s_{j}} \over 1-e^{-\alpha_{\sigma(j)}\bar s}}\int_{S^k}\prod_{j=1}^k
m(dy_{j})f_{\sigma(j)}(y_{j})
 \Gamma_{y_{j-1},y_{j}}(ds_{j})
 \nonumber\\
=&&\int_{s_{1}>0,\ldots,s_{k}>0}g(\bar s)
  e^{-\alpha_{\sigma'(1)}\bar s} \nn\\
&&\hspace{.2 in} \prod_{j=1}^k {e^{-\alpha_{\sigma'(j)}\hat s_{j}} \over 1-e^{-\alpha_{\sigma'(j)}\bar s}}
\int_{S^k}\prod_{j=1}^k m(dy_{j})f_{\sigma'(j)}(y_{j})
 \Gamma_{y_{j-1},y_{j}}(ds_{j}) \nn
\end{eqnarray}
whenever $(\sigma'(1),\sigma'(2),\ldots, \sigma'(k))=(\sigma(k),\sigma(1),\ldots, \sigma(k-1))$.

Note that 
\begin{eqnarray}
&&M_{\si'}(ds_{1},\ldots,ds_{k}):=\int_{S^k} \prod_{j=1}^k  m(dy_{j})f_{\sigma'(j)}(y_{j})
 \Gamma_{y_{j-1},y_{j}}(ds_{j})
\label{retr}\\
&&= \int_{S^k} f_{\sigma'(1)}(y_{1})
 \Gamma_{y_{k},y_{1}}(ds_{1})f_{\sigma'(2)}(y_{2})
 \Gamma_{y_{1},y_{2}}(ds_{2})\cdots \nonumber\\
 &&\hspace{2 in}\cdots f_{\sigma'(k)}(y_{k})
 \Gamma_{y_{k-1},y_{k}}(ds_{k}) \prod_{j=1}^k  m(dy_{j})    \nonumber\\
&&= \int_{S^k}  f_{\sigma(k)}(y_{1})
 \Gamma_{y_{k},y_{1}}(ds_{1})f_{\sigma(1)}(y_{2})
 \Gamma_{y_{1},y_{2}}(ds_{2})\cdots \nonumber\\
 &&\hspace{2 in} \cdots f_{\sigma(k-1)}(y_{k})
 \Gamma_{y_{k-1},y_{k}}(ds_{k})  \prod_{j=1}^k  m(dy_{j}),    \nonumber
 \eea
 and relabeling the $y_{j}$'s this is
 \bea
&&= \int_{S^k} f_{\sigma(k)}(y_{k})
 \Gamma_{y_{k-1},y_{k}}(ds_{1})f_{\sigma(1)}(y_{1})
 \Gamma_{y_{k},y_{1}}(ds_{2})\cdots \label{retrc}\\
 &&\hspace{1.5 in} \cdots f_{\sigma(k-1)}(y_{k-1})
 \Gamma_{y_{k-2},y_{k-1}}(ds_{k}) \prod_{j=1}^k  m(dy_{j})    \nonumber\\
 &&=\int_{S^k} \prod_{j=1}^k  m(dy_{j})f_{\sigma(j)}(y_{j})
 \Gamma_{y_{j-1},y_{j}}(ds_{j+1})= M_{\si}(ds_{2},\ldots,ds_{k},ds_{1}),\nn
\end{eqnarray}
where $s_{k+1}=s_{1}$.
Furthermore, (recall that $\hat s_{1}=0$), 
\begin{eqnarray}
&&\alpha_{\sigma'(1)}\bar s+\sum_{j=2}^k\alpha_{\sigma'(j)}\hat s_{j}
\label{retrx}\\
&& =\alpha_{\sigma(k)}\bar s+\sum_{j=2}^k\alpha_{\sigma(j-1)}\(\sum_{l=2}^{j}s_{l}\)  \nonumber\\
&& =\alpha_{\sigma(k)}\bar s+\sum_{i=1}^{k-1}\alpha_{\sigma(i)}\(\sum_{l=2}^{i+1}s_{l}\)  \nonumber\\
&& =\alpha_{\sigma(k)}(s_1+\cdots+s_k)+\sum_{i=1}^{k-1}\alpha_{\sigma(i)}\(\sum_{l=1}^{i}s_{l+1}\)=\sum_{i=1}^{k}\alpha_{\sigma(i)}\(\sum_{l=1}^{i}s_{l+1}\).  \nonumber
\end{eqnarray}
But also, using $\sum_{j=1}^{k}\al_{j}=\sum_{j=1}^{k}\al_{\si (j)}$,
 \be
(\alpha_1+\cdots\alpha_k)s_{1}+\sum_{j=2}^k \alpha_{\sigma(j)}\hat s_{j}=\sum_{j=1}^k \alpha_{\sigma(j)}\(\sum_{l=1}^{j}s_{l}\)\label{p10.1n}
\end{equation}
where $s_{k+1}=s_{1}$.
Combining (\ref{retr})-(\ref{p10.1n}), we see that (\ref{p8m}) is the claim that  
\begin{eqnarray}
&&\int_{s_{1}>0,\ldots,s_{k}>0} g(\bar s)
\exp \(- \sum_{j=1}^k \alpha_{\sigma(j)}\(\sum_{l=1}^{j}s_{l}\)\)\label{p8mu}\\
&&\hspace{2 in} \prod_{j=1}^k {1 \over 1-e^{-\alpha_{\sigma(j)}\bar s}}M_{\si}(ds_{1},\ldots,ds_{k})
 \nonumber\\
=&&\int_{s_{1}>0,\ldots,s_{k}>0}g(\bar s)
\exp \(- \sum_{j=1}^k \alpha_{\sigma(j)}\(\sum_{l=1}^{j}s_{l+1}\)\) \nn\\
&&\hspace{2 in} \prod_{j=1}^k {1 \over 1-e^{-\alpha_{\sigma'(j)}\bar s}} M_{\si}(ds_{2},\ldots,ds_{k},ds_{1}) \nn
\end{eqnarray}
where $s_{k+1}=s_{1}$, and this claim follows immediately from the relabeling $(s_{1},\ldots,s_{k})\to (s_{2},\ldots,s_{k},s_{1})$.
This establishes (\ref{p8m}) and hence (\ref{pj1}). \qed

 \section{The restriction property}\label{sec-rest}

Let $B\subseteq S$ be open and set
\begin{equation}
T_{B^{c}}=\inf \{t\geq 0\,|\ X_{t}\in  B^{c} \}.\label{6.50}
\end{equation}
  Let
         \begin{equation}
\wt X_t(\om)=\left\{\begin{array}{ll} X_t(\om) & \mbox{if
$t<T_{B^{c}}$}\\
         \Delta & \mbox{otherwise}.
         \end{array}\right.\label{mpfr.20}
         \end{equation}
        Clearly, $t\mapsto \wt X_t$ is  right continuous. With
\begin{equation}
\wt P_tf(x)= E^x(f( X_t)1_{\{t<T_{B^{c}}\}}),\label{mpfr.21}
\end{equation} and
  \begin{equation}
\wt \th_t(\om)=\left\{\begin{array}{ll} \th_t(\om) &
\mbox{if $t<T_{B^{c}}( \om)$}\\
         \Delta & \mbox{otherwise,}
         \end{array}\right.\label{mpfr.20a}
         \end{equation} 
we show in \cite[Section 4.5]{book} that $\wt X=(\Om,
\, \mathcal{G}_t, \mathcal{G},\wt  X_t,\wt  \th_t,\wt  P^x
)$ is a Borel right process with state space $B$ and potential densities
  \be
         \wt{u}  (x,y)= u  (x,y)-E^x\( \, u  \(
X_{T_{B^{c}} },y\)\), \hspace{ .4in}x,y\in B,\label{mpfr.27}
          \ee 
with respect to the measure $m(dx)$ restricted to $B$. Here we have used the convention that $u(\Delta, y)=0$ and that $X_t(\omega) = \Delta$ when $t=+\infty$. It follows as before that uniformly in $x$, $ \wt{u}  (x,y)$ is locally bounded and continuous in $y$. 

Let $\{L^{x}_{t},\,(x,t)\in S\times R_{+}\}$ be the  family of local times for $X$ used in the construction of $\mu$. Set $\wt L^{x}_{t}=L^{x}_{t\wedge T_{B^{c}}}$ for $x\in B$. It is easy to see that $\wt L^{x}_{t}$ is a CAF for $\wt X$ and
\begin{eqnarray}
\wt E^{x}\(\wt L^{y}_{\ff}\)&=&E^{x}\(L^{y}_{T_{B^{c}}}\)
\nn\\
&=&    E^{x}\(L^{y}_{ \ff }\)-E^{x}\(L^{y}_{\ff}\circ\th_{T_{B^{c}}}\)\nonumber\\
&=&    u  (x,y)-E^x\( \, u  \(
X_{T_{B^{c}}},y\)\)= \wt{u}  (x,y).\nonumber
\end{eqnarray}
It follows that $\{\wt L^{x}_{t},\,(x,t)\in B\times R_{+}\}$ are local times for $\wt X$.
We can then define the loop measure  $\wt \mu$  for $\wt X$ by the formula 
\begin{equation}
\int F\,d\wt \mu=\int_{B}\wt P^{x}\(\int_{0}^{\ff}{1\over t}\,F\circ k_{t}\, \,d_{t}\wt L_{t}^{x}\)\,dm (x).  \label{ls.335}
\end{equation}
(In our definition (\ref{ls.3}) of $\mu$ we assumed that $X$ had continuous potential densities. We do not know if $  \wt{u}  (x,y)$ is continuous in $x$. However, the continuity of $  u(x,y)$ was only used to guarantee a nice family of local times for $X$, and by the above this is inherited by 
$\{\wt L^{x}_{t},\,(x,t)\in B\times R_{+}\}$).  
\bt[The Restriction Property]\label{theo-rest}
\begin{equation}
\mu(F; T_{B^c}=\infty)=\wt\mu(F).\label{6.51}
\end{equation}
\et

Note that $B^{c}=S-B$ does not contain $\De$.

{\bf  Proof of Theorem \ref{theo-rest}: } It suffices to prove this for $F$ of the form 
$\prod_{j=1}^{k}f_{j}(X_{t_{j}})$ with $t_{1}<\cdots<t_{k}$. 
Since $f_{k}(\De)=0$,  and $1_{\{T_{B^c}=\infty\}}\circ k_{t}=1_{\{t\leq T_{B^c}\}}$ we have
\begin{equation}
1_{\{T_{B^c}=\infty\}}\circ k_{t}\,\,\prod_{j=1}^{k}f_{j}(X_{t_{j}})\circ k_{t}=1_{\{t_{k}<t\leq T_{B^c}\}}\prod_{j=1}^{k}f_{j}(X_{t_{j}}).\label{6.51r}
\end{equation}
Hence
\bea
&&
\mu\(   \prod_{j=1}^{k}f_{j}(X_{t_{j}}); T_{B^c}=\infty\)\label{6.51s}\\
&&
=\int_{S}P^{x}\(\int_{t_{k}}^{T_{B^c}}{1\over t}\,\prod_{j=1}^{k}f_{j}(X_{t_{j}})\,dL_{t}^{x}\)\,dm (x) \nn\\
&&
=\int_{B}P^{x}\(\prod_{j=1}^{k}f_{j}(X_{t_{j}})1_{\{t_{k}<T_{B^c}\}}\,\int_{t_{k}}^{\ff}{1\over t}\,d\wt L_{t}^{x}\)\,dm (x) \nn\\
&&
=\int_{B}P^{x}\(\prod_{j=1}^{k}f_{j}(X_{t_{j}})1_{\{t_{k}<T_{B^c}\}}\,\(\int_{t_{k}}^{\ff}{1\over t}\,d_{t}\wt L_{t-t_{k}}^{x}\)\circ \th_{t_{k}}\)\,dm (x) \nn\\
&&
=\int_{B}P^{x}\(\prod_{j=1}^{k}f_{j}(X_{t_{j}})1_{\{t_{k}<T_{B^c}\}}\,E^{X_{t_{k}}}\(\int_{t_{k}}^{\ff}{1\over t}\,d_{t}\wt L_{t-t_{k}}^{x}\)\)\,dm (x). \nn
\eea
 But this is clearly
 \begin{equation}
\int_{B}\wt P^{x}\(\prod_{j=1}^{k}f_{j}(\wt X_{t_{j}})\,\wt E^{\wt X_{t_{k}}}\(\int_{t_{k}}^{\ff}{1\over t}\,d_{t}\wt L_{t-t_{k}}^{x}\)\)\,dm (x) \label{6.51t}
 \end{equation}
which is precisely what we obtain from $\wt\mu\( \prod_{j=1}^{k}f_{j}(X_{t_{j}})\)$
by proceeding as in (\ref{6.51s}).
\qed

\section{Transformations of the loop measure}\label{sec-trans}

\subsection{Mappings of the state space}\label{sec-spacetrans}

Let $\bar S$  be another locally compact 
 topological space
 with a countable base, and let $f: S\mapsto \bar S$ be a topological isomorphism. Then
 \begin{equation}
\bar P_{t}(x,g)=P_{t}(f^{-1}(x),g\circ f).\label{6.70a}
\end{equation}
forms a Borel transition semigroup on $\bar S$.
 Let   $\bar \Om$ be the set of right continuous paths $\om$ in $\bar S_{\De}=\bar S\cup \De$ with  $\De\notin \bar S$, and   such that   $\om(t)=\De$ for all $t\geq \ze=\inf \{t>0\,|\om(t)=\De\}$.   Then with $\bar X_{t}(\om)=\om(t)$  it follows from \cite[Section 13]{S}  that  $\bar X\!=\!
(\bar \Om, \bar  \FF_{t},\bar  X_t,\th_{t} ,\bar P^x
)$ is a Borel right process.
Furthermore,
\bea
\bar U(x,g)&=&\int_{0}^{\ff}\bar P_{t}(x,g)\,dt=\int_{0}^{\ff}P_{t}(f^{-1}(x),g\circ f)\,dt\label{6.70}\\
&=&\int_{S}u(f^{-1}(x),y)g\circ f(y)\,dm(y)\nn\\
&=&\int_{\bar S}u(f^{-1}(x),f^{-1}(z))g(z)\,df_{\ast}m(z).\nn
\eea
Thus $\bar X$ has continuous potential densities $\bar u(x,y)=u(f^{-1}(x),f^{-1}(y))$ with respect to the measure $f_{\ast}m$.

If we let $\bar f:\Om\mapsto \bar \Om$ be defined as $\bar f(\om)(t)=f(\om(t))$, it follows that  
\begin{equation}
\bar P^x \(F\)=P^{f^{-1}(x)} \(F\circ \bar f\,\,\).\label{6.70b}
\end{equation}
Note further that $\bar L^{y}_{t}=L^{f^{-1}(y)}_{t}\circ\bar  f^{-1}$ is a CAF for 
$\bar  X$ with
\begin{equation}
\bar P^x \(\bar L^{y}_{\ff}\)=P^{f^{-1}(x)} \(L^{f^{-1}(y)}_{\ff}\)=u(f^{-1}(x),f^{-1}(y))=\bar u(x,y),\label{6.71}
\end{equation}
so that $\{\bar L^{y}_{t},\,(y,t)\in S'\times R_{+}\}$ are local times for $\bar X$. Let $\bar\mu$ be the loop measure for $\bar X$.

\bt\label{theo-spat}
\begin{equation}
\bar  f_{\ast}\mu\(F\)=\bar  \mu\(F\). \label{6.75}
\end{equation}
\et
{\bf  Proof of Theorem \ref{theo-spat}: }
\bea
\bar  \mu\(F\)&=&\int_{S'}\bar P^{x}\(\int_{0}^{\ff}{1\over t}\,F\circ k_{t}\, \,d_{t}\bar L_{t}^{x}\)\,df_{\ast}m (x)  \label{6.72}\\
&=& \int_{S'}  P^{f^{-1}(x)}\(\int_{0}^{\ff}{1\over t}\,F\circ k_{t}\circ\bar  f\, \,d_{t}  L_{t}^{f^{-1}(x)}\)\,df_{\ast}m (x)\nn\\
&=& \int_{S}  P^{x}\(\int_{0}^{\ff}{1\over t}\,F\circ k_{t}\circ\bar  f\, \,d_{t}  L_{t}^{x}\)\,dm (x)\nn\\
&=&\mu\(  F\circ \bar f \)=\bar  f_{\ast}\mu\(F\).\nn
\eea
\qed

\subsection{Unit Weights}

We say that a random variable $T\geq 0$ is a unit weight if
\begin{equation}
\int_{0}^{\ze}T\circ \rho_{u}\,du=1.\label{6.40}
\end{equation}
Of course, since $\ze$   is invariant under  loop  rotation,  $1/\ze$ is an example of a unit weight.
(\ref{c6.84}) will provide another example,  which is be used   in the proof of Theorem \ref{theo-cconfinv} to determine how the loop measure transforms under a time change.

 Let $\mathcal{I}_{\rho}(X)$ be the collection of measurable functions on $\Om $ which are invariant under   loop  rotation.   
\bl If $T$ is a unit weight then  for all $F\in \mathcal{I}_{\rho}(X)$
\begin{equation}
\mu(F)= \int_{S}\,Q ^{x,x}\(T\,F  \)\,dm (x).\label{6.42} 
\end{equation}
\el
{\bf  Proof: }
By invariance of $\mu$ we have that for each $u>0$ and $F\in \mathcal{I}_{\rho}(X)$
\begin{equation}
\mu (T\circ \rho_{u}\,\,F)=\mu (T\,F).\label{6.42a}
\end{equation}
Since $\ze$   is invariant under  loop  rotation, this implies that for any $u>0$
\begin{equation}
\mu (T\circ \rho_{u}\,1_{\{u<\ze\}}\,F)=\mu (T\,1_{\{u<\ze\}}\,F).\label{6.42b}
\end{equation}
Hence
\begin{equation}
\int_{0}^{\ff}\mu (T\circ \rho_{u}\,1_{\{u<\ze\}}\,F)\,du=\int_{0}^{\ff}\mu (T\,1_{\{u<\ze\}}\,F)\,du.\label{6.43}
\end{equation}
This shows that
\begin{equation}
\mu \(\int_{0}^{\ze}T\circ \rho_{u}\,du\,F\)= \mu \(T\,F\,\int_{0}^{\ff}1_{\{u<\ze\}}\,du\)
= \mu \(T\,F\,\ze\).\label{6.43a}
\end{equation}
Using (\ref{6.40}) and (\ref{6.41}) our Lemma follows.\qed

\subsection{Time change by the inverse of a CAF}\label{sec-time}

Consider
\begin{equation}
A_{t}=\int_{S} L^{x}_{t}\,d\nu_{A}(x)\label{101.2}
\end{equation}
where $\nu_{A}$ is a Borel measure on $S$. We suppose that  $P^x(A_t=\infty, t<\zeta)=0$ for all $x\in S$ and $t>0$. (This is the case, for instance, if $\nu_{A}(K)<\infty$ for each compact $K\subset S$.)  By the argument at the beginning of the proof of Lemma \ref{lem-ltmom},  (\ref{101.2}) defines a CAF of $X$. 
Let $S_{A}$ denote the fine support of $A$; that is, the set of $x\in S$ such that $P^{x}(R=0)=1$ where $R=\inf \{t>0\,|\, A_{t}>0\}$, see \cite[Section 64]{S}. Because $m$ is a reference measure and $v(x):=E^x(\exp(-R))$ is a $1$-excessive function,  $S_A =\{x\in S: v(x)=1\}$ is a Borel subset of $S$; see \cite[Prop.~V(1.4)]{BG}.

Let $\tau_{t}$ be the right continuous inverse of $A_{t}$, and let $Y_{t}=X_{\tau_{t}}$. Then  $Y\!=\!
(\Om,  \GG_{t}, Y_t,\wh \th_{t} ,P^x
)$ is a Borel right process  with state space $S_{A}$ and lifetime $A_{\ze}$, see \cite[Theorem 65.9]{S} for details, noting that \cite[(60.4)]{S} applied to $\exp(-A_t)$ allows us to assume that $A$ is a perfect CAF. Here $\wh \th_{t}(\om)=\th_{\tau_{t}(\om)}(\om)$. 
 Using the change of variables formula, \cite[Chapter 6, (55.1)]{DM2}, we see that 
\bea
E^{x}\(\int_{0}^{\ff} f\(Y_{t}\)\,dt\)&=&E^{x}\(\int_{0}^{\ff} f\(X_{\tau_{t}}\)\,dt\)\\\label{6.82adn}&=&
E^{x}\(\int_{0}^{A_{\ze}} f\(X_{\tau_{t}}\)\,dt\)\nn\\
&=&E^{x}\(\int_{0}^{\ff} f\(X_{s}\)\,dA_{s}\)\nn\\
&=& \int u(x,y)f(y)\,\nu_A(dy),\nn
\eea
so that $Y_{t}$ has continuous potential densities $u(x,y)$ with respect to the measure $\nu_A(dy)$ on $S_{A} $. (In the last step we used the fact that for any measurable function $h_{s}$, we have 
$\int_{0}^{\ff} h_{s}\,dA_{s}=\int \(\int_{0}^{\ff} h_{s}\,dL^{y}_{s}\)\,\nu_A(dy)$. It suffices to verify this for functions of the form  $h_{s}=1_{[0,t]}(s)$, in which case it is obvious). 
Furthermore, since $S_{A}$ is the fine support of $A$, $L^{x}_{\tau_{t}}$ is continuous in $t$ for each $x\in S_{A}$, see \cite[p. 1659]{G}, and of course $E^{y}\(L^{x}_{\tau_{\ff}}\)=u(y,x)$. It follows that  $\{L^{x}_{\tau_{t}},\,(x,t)\in S_{A}\times R_{+}\}$ is the  family of local times for $Y$. See \cite[Theorem 65.6]{S} for additivity. 

It will be convenient to use the canonical notation $ X'\!=\!
( \Om,   \FF'_{t},  X'_t,\th_{t} ,P'^x)$ for $Y$. Thus $X'_{t}(\om)=\om(t)$, which is the same as $X_{t}(\om)$, but we use the notation $X'_{t}$ to emphasize that it is associated with the measures $P'^x$ which we now define. If we set $g(\om)(t)=\om(\tau_{t}(\om))$ we have $Y_{t}=X_t\circ g$ and put
\begin{equation}
P'^x\(F\)=P^x\(F\circ g\).\label{con.1}
\end{equation}
%By the above $\{L'^{x}_t:=L^{x}_{t}\circ g,\,(x,t)\in S_{A}\times R_{+}\}$ is the  family of local times %for $X'$.   
Using   \cite[(62.20)]{S}, compare (\ref{nit4.4}), we see that if $t_{1}<\cdots<t_{n}$,
\begin{eqnarray}
 g_{\ast}Q^{x,x}\(\prod_{j=1}^{n}f_{j}\(X_{t_{j}}\)\)
&=& Q^{x,x}\(\prod_{j=1}^{n}f_{j}\(X_{\tau_{t_{j}}}\)\)  
\label{c6.83a}\\
&=& P^{x}\(\prod_{j=1}^{n}f_{j}\(X_{\tau_{t_{j}}}\)u(X_{\tau_{t_{n}}},x)\)  \nonumber\\
&=& P'^{x}\(\prod_{j=1}^{n}f_{j}\(X'_{t_{j}}\)u(X'_{t_{n}},x)\).  \nonumber
\end{eqnarray}
Let  $Q'^{x,x}$ be the measure in (\ref{ls.3n}) associated with $X'$. Using (\ref{nit4.4}) and the fact that $X'$ also has potential densities $u(x,y)$ we have that if $x\in S_{A}$
\begin{eqnarray}
&&P'^{x}\(\prod_{j=1}^{n}f_{j}\(X'_{t_{j}}\)u(X'_{t_{n}},x)\)= Q'^{x,x}\(\prod_{j=1}^{n}f_{j}\(X'_{t_{j}}\)\). \nonumber
\end{eqnarray} 
Hence for all  measurable $F$
\begin{equation}
g_{\ast}Q^{x,x}\(F\)=Q'^{x,x}\(F\),\hspace{.2 in}\forall x\in S_{A}.\label{c6.83}
\end{equation}

Before considering general $\nu_{A}$'s, we first study  the special case where the measure $\nu_{A}$ is equivalent to $m$. Thus 
$\nu_{A}(dx)=h(x)m(dx)$ where  $h$ is a measurable  function on $S$ with $0<h(x)<\ff$ for all $x$.  It follows from (\ref{odf}) that
\begin{equation}
A_{t}=\int_{0}^{t}h\(X_s\)\,ds,\label{c6.82}
\end{equation}
and thus  $S_{A}=S$.
Let $\mu, \mu'$ be the loop measures for $X,X'$ 
respectively.
\bt\label{theo-cconfinv} If $\nu_{A}(dx)=h(x)m(dx)$ where  $h$ is a measurable  function on $S$ with $0<h(x)<\ff$ for all $x$, then
\begin{equation}
g_{\ast}\mu\(F\)=\mu'\(F\), \hspace{.2 in}\forall F\in\mathcal{I}_{\rho}\(X'\). \label{c6.800}
\end{equation}
\et 

{\bf  Proof of Theorem \ref{theo-cconfinv}: }
Define the unit weight
\begin{equation}
T(\om)={h(\om (0)) \over A_{\zeta}(\om)}.\label{c6.84}
\end{equation}
By (\ref{6.42}) we have $\mu\(F\)= \int_{S}   Q^{x,x}\(T\,\,F \)\,m(dx) $ for all  $F\in\mathcal{I}_{\rho}\( X\)$. It is easy to see that $F\in\mathcal{I}_{\rho}\( X'\)$ implies that $F\circ g\in\mathcal{I}_{\rho}\( X\)$.
Noting  that $A_{\zeta}=\ze\circ g$, and using (\ref{c6.83})
\begin{eqnarray}
g_{\ast}\mu\(F\)&=&\mu\(F\circ g\)
\label{c6.85}\\
&=& \int_{S}   Q^{x,x}\(T\,\,F\circ g \)\,m (dx)  \nonumber\\
&=& \int_{S}   \,Q^{x,x}\( {1 \over \ze\circ g }F\circ g \)h(x)\,m (dx)  \nonumber\\
&=& \int_{S}  Q'^{x,x}\( {1 \over \ze} \,F  \)h(x)\,m (dx)  =\mu'\(F\).\nonumber
\end{eqnarray}
The last equality used (\ref{6.41}) and the fact that $\nu_{A}(dx)=h(x)m(dx)$.
\qed

We next show how to combine Theorems  \ref{theo-spat} and \ref{theo-cconfinv}.
Let $S'$  be another locally compact 
 topological space
 with a countable base, and let $f: S\mapsto S'$ be a topological isomorphism.  With $h$ as above, let  $m_{S'}$ be the measure  on $S' $ defined by 
\begin{equation}
m_{S'}(dy):=f_{\ast}\(h\,\, m_{S} \)(dy).\label{6.81}
\end{equation} 
%Note that  
%\begin{equation}
%h(x)m_{S}(dx)=(f^{-1})_{\ast} m_{S'}(dx).\label{6.81f}
%\end{equation}
It follows from the discussion in sub-section \ref{sec-spacetrans} and the present sub-section that if we set
$\bar X'_{t}:=f\(X_{\tau_{t}}\)=f\(Y_{t}\)$ and $\{\bar P'^x, x\in S'\}$ the measures induced by $\{P^x, x\in S\}$, then $\bar X'=\!
(\bar  \Om, \bar  \FF_{t}, \bar X'_t,\th_{t} , \bar  P'^x
)$ is a Borel right process with continuous potential densities
\begin{equation}
\bar u'(x,y)=u(f^{-1}(x),f^{-1}(y))\label{6.81g}
\end{equation}
  with respect to the measure $f_{\ast}\(h\,\, m_{S} \)=m_{S'}$ on $S'$.

Set $f^{\sharp}(\om)(t)=f(\om(\tau_{t}))$ and let $\mu, \bar  \mu'$ be the loop measures for $X, \bar X'$ 
respectively. Combining Theorems  \ref{theo-spat} and \ref{theo-cconfinv} we obtain
\begin{corollary}\label{cor-confinv}
\begin{equation}
f^{\sharp}_{\ast}\mu\(F\)=\bar \mu'\(F\), \hspace{.2 in}\forall F\in\mathcal{I}_{\rho}\(X'\). \label{6.800}
\end{equation}
\end{corollary}
\begin{remark} {\rm    Let $D,D' $ be  two simply connected domains in the complex plane and let $f$ be a conformal map from $D$ onto $D'$. Let $X$ be Brownian motion in D. Since the potential density for $X$ with respect to $\la_{D}$, Lebesgue measure on $D$, is not continuous, (it has a logarithmic singularity on the diagonal), $X$ does not fit into the framework of this paper. Nevertheless, we argue by analogy. Let $h(x)=|f'(x)|^{2}$. Then $\bar X'$ is a Brownian motion in $D'$, and $ f_{\ast}\,(h\,\la_{D})(dy)=\la_{D'}(dy)$. It follows formally  that (\ref{6.800}) would yield \cite[Proposition 5.27]{Lawler}, the conformal invariance of Brownian loop measures.   }
 \end{remark}

We now turn to a general CAF as in (\ref{101.2}), Our results are not as complete as (\ref{c6.800}), but see the Remark following the proof of Theorem \ref{theo-gtc}.

For any $B\subseteq S$ let $\mathcal{L}_{B}(X)$ be the $\si$-algebra generated by the total  local times $\{L^{x}_{\ff},\,x\in B\}$ of $X$, and  let $\mu'$ be the loop measure for $X'$.

  \bt\label{theo-gtc}
\begin{equation}
g_{\ast}\mu\(F\)=\mu'\(F\), \hspace{.2 in}\forall F\in\mathcal{L}_{S_{A}}\(X'\). \label{101.0}
\end{equation}
\et

{\bf  Proof of Theorem \ref{theo-gtc}: }By   Lemma 2.2
\be
\mu(L^x_\infty\cdot F) =Q^{x,x}(F),\hspace{.2 in}\forall  F\in \mathcal{L}_{S}(X),\,\, x\in S.  \label{101.1}
\ee 
Recall that $L'^{x}_\ff=L^{x}_\ff\circ g$ so that
\begin{equation}
\mathcal{L}_{S_{A}}\(X'\)=\mathcal{L}_{S_{A}}\(X\)\circ g.\label{sh.9}
\end{equation}
Consider $F\in\mathcal{L}_{S_{A}}\(X\)$. Since $A_{\ze}\in \mathcal{L}_{S_{A}}(X)$ and   $A_{\ze}>0$, $P_{x}$ a.s. for all $x\in S_{A}$, by replacing $F$ in (\ref{101.1}) by $F/A_{\ze}$ and then integrating with respect to  $d\nu_{A}(x)$ we can deduce immediately that
\be
\mu(F) =\int_{S_{A}} Q^{x,x}\({F \over A_\zeta}\)\nu_A(dx), \hspace{.2 in}\forall F\in\mathcal{L}_{S_{A}}\(X\).\label{101.3}
\ee  

 Although $S_{A}$ may not be locally compact, $X'$ inherits from $X$ all the   properties needed to define $ \mu' $ as in (\ref{ls.3}),
and it then follows as in  (\ref{ls.3n})  that
\begin{equation}
\mu'(F)= \int_{S_{A}}\,Q'^{x,x}\( {1 \over  \zeta}\,F \)\,d\nu_A (x).\label{6.41n}
\end{equation} 
By (\ref{c6.83})  this shows  that 
\begin{equation}
\mu'(F)= \int_{S_{A}}\,Q^{x,x}\( {1 \over \zeta\circ g}\,F\circ g \)\,d\nu_A (x).\label{6.41n1}
\end{equation} 
Noting  that $A_{\zeta}=\ze\circ g$,   (\ref{101.3}) and (\ref{sh.9}) then imply our Theorem.
 \qed

 \begin{remark}{\rm  For $x_{1},\ldots, x_{n}\in S$ we define the multiple local time
\bea
&&
L_{t}^{x_{1},\ldots, x_{n}}\label{6.90}\\
&&=\sum_{j=1}^{n} \int_{0\leq s_{1}\leq\ldots\leq s_{n}\leq t}dL^{x_{j}}_{s_{1}}dL^{x_{j+1}}_{s_{2}}\cdots
dL^{x_{n}}_{s_{n-j+1}}dL^{x_{1}}_{s_{n-j+2}}dL^{x_{2}}_{s_{n-j+3}}\cdots dL^{x_{j-1}}_{s_{n}},\nn
\eea
that is, we measure $n$-tuples of times   in which $x_{1},\ldots, x_{n}$ are visited in cyclic order. If $n=2$ and $x_{1}\not=x_{2}$, then $L_{t}^{x_{1},  x_{2}}=L_{t}^{x_{1}}L_{t}^{ x_{2}}$, but in general $L_{t}^{x_{1},\ldots, x_{n}}$ is not a product of the corresponding local times. Let $\mathcal{M}(X)$ denote the $\si$-algebra generated by the multiple local times.  When   $\mbox{Supp}\,(\nu_A) =S $ we can show that   (\ref{101.0}) holds for all  $F\in\mathcal{M}(X)=\mathcal{M}(X')$. When $S$ is finite, it is known that $\mathcal{M}(X)=\mathcal{I}_{\rho}\(X\)$, \cite[p. 24]{Le Jan1}.   For diffusions, see \cite{Lupu}, especially Corollary 2.9, and for more general processes see \cite{Chang}.}
\end{remark}
 
We leave to the interested reader the task of combining Theorem \ref{theo-gtc} with spatial transformations as in Corollary \ref{cor-confinv}.

\def\noopsort#1{} \def\printfirst#1#2{#1}
\def\singleletter#1{#1}
            \def\switchargs#1#2{#2#1}
\def\bibsameauth{\leavevmode\vrule height .1ex
            depth 0pt width 2.3em\relax\,}
\makeatletter
\renewcommand{\@biblabel}[1]{\hfill#1.}\makeatother
\newcommand{\bysame}{\leavevmode\hbox to3em{\hrulefill}\,}

 \def\wh{\widehat}
\def\ol{\overline}

\begin{thebibliography}{10}

\bibitem{BG}
R. M. Blumenthal, and R. K. Getoor, (1968)
\newblock{\em Markov Processes and Potential Theory}, 
\newblock Academic Press, New York.

\bibitem{Bertoin} J.  Bertoin, {\em Random Fragmentation and Coagulation Processes}, Cambridge University Press, New
York,  (2006).

 \bibitem {Chang} Y. Chang,
Multi-occupation field generates the Borel sigma-field of loops. arXiv:1309.1558

\bibitem {DM}
C. Dellacherie,  and P.-A. Meyer,  (1978).
\newblock {\em Probabilities et Potential}.
\newblock Paris: Hermann.

\bibitem {DM2}
C. Dellacherie,  and P.-A. Meyer,  (1982).
\newblock {\em Probabilities and Potential B}.
\newblock North Holland Publishing Company, Amsterdam.


\bibitem {D80}
E. B. Dynkin,
\newblock Minimal excessive functions and measures.  \newblock {\em Trans.\ Amer.\ Math.\ Soc.} {\bf 258}  (1980), 217--244. 

\bibitem {D83}
E. B. Dynkin,
\newblock Local times and quantum fields, In   {\em Seminar on Stochastic Processes, 1983}, pp.~64--84.  Progress in Probability,  Birkh\"auser, Boston (1984).  

\bibitem {D84}
E. B. Dynkin,
\newblock{G}aussian and non--{G}aussian random fields associated with {M}arkov
    processes.    {\em J.\  Funct.\ Anal.} {\bf 55}  (1984), 344--376.   




\bibitem {EK}
N. Eisenbaum  and  H. Kaspi,
\newblock On permanental processes, {\em Stochastic Processes and their Applications},
{\bf 119}  (2009),  1401--1415.

\bibitem {FPY}
P. Fitzsimmons, J. Pitman,  and M. Yor,  
\newblock Markovian bridges: construction, Palm interpretation, and splicing.  
{\em Seminar on Stochastic Processes}, 1992, pp.~101--134. Birkh\"auser, Boston (1993).

\bibitem {G}
R. K. Getoor,
\newblock Some remarks on continuous additive functionals,  {\em Ann. Math. Statist.},
{\bf 38}  (1967),  1655--1660.



\bibitem{Le Jan} Y. Le Jan, \newblock
  Markov loops and renormalization, {\em Ann. Probab.},
{\bf 38}  (2010),  1280--1319. 

\bibitem{Le Jan1} Y. Le Jan, \newblock
{\em Markov paths, loops and fields. }   \'{E}cole d'\'{E}t\'{e} de Probabilit\'{e}s de Saint-Flour XXXVIII - 2008. Lecture Notes in Mathematics 2026. 
Springer-Verlag, Berlin-Heidelberg, (2011).

\bibitem{LMR1} Y. Le Jan, M. B.  Marcus and J.~Rosen,
\newblock  Permanental fields,  loop soups and continuous additive functionals, {\em Ann. Probab.}, to appear.
\newblock   arxiv.org/pdf/1209.1804.pdf 

\bibitem{LMR2} Y. Le Jan, M. B.  Marcus and J.~Rosen,
\newblock    Intersection local times, loop soups and permanental Wick powers.
\newblock arxiv.org/pdf/1308.2701.pdf



 

\bibitem{K} J. F. C. Kingman, {\em Poisson Processes}, Oxford Studies in Probability, Clarendon Press, Oxford, (2002).



\bibitem{LL} G. Lawler and V. Limic, {\em Random Walk: A Modern Introduction}, Cambridge University Press, New
York,  (2009).

   \bibitem{Lawler} G. Lawler, {\em Conformally Invariant
Processes in the Plane}, Math. Surveys and Monographs, v. 114,  AMS, (2005).


\bibitem{LF}
G. Lawler and J. Trujillo Ferreras, \newblock Random walk loop soup, \newblock {\em Trans.\ Amer.\ Math.\ Soc. } {\bf 359}  (2007), 565--588.

\bibitem{LW}
G. Lawler and W. Werner,  \newblock The Brownian loop soup, \newblock {\em Prob.\ Theory.\ Rel.\ Fields} {\bf 44}  (2004), 197--217.

\bibitem{Lupu} T. Lupu,
\newblock  Poissonian ensembles of loops of one-dimensional diffusions.
\newblock   arxiv.org/pdf/1302.3773.pdf
 
\bibitem{book} M. B.  Marcus and J.~Rosen, {\em Markov Processes,
Gaussian Processes and Local Times}, Cambridge University Press, New
York,  (2006).



  
   \bibitem{S} M. Sharpe, {\em General theory of Markov
processes},  Acad. Press, New York, (1988).


\bibitem {Sy}
K. Symanzik,
\newblock Euclidean quantum field theory, In   {\em Local Quantum Theory}, (R. Jost, ed.). pp.~152--226.  Acad. Press, New York, (1967).  


\bibitem{Sz1} A.-S. Sznitman,
\newblock  Topics in occupation times and Gaussian free fields.
 {\it Zurich Lectures in Advanced Mathematics}, EMS, Zurich, 2012.
 


\bibitem{VJ}  D.  Vere-Jones,
  \newblock Alpha-permanents and their applications to multivariate gamma, negative binomial and ordinary binomial distributions.
    {\em New Zealand J.  Math.}, {\bf 26}  (1997), 125--149.



\end{thebibliography}
\end{document}